\documentclass[a4paper]{amsart}
\usepackage[T1]{fontenc}
\usepackage[utf8]{inputenc}
\usepackage[british]{babel}
\usepackage{amsfonts,amsmath,amssymb,amsthm,mathtools,esint}
\usepackage{eucal} 
\usepackage[dvipsnames]{xcolor}
\usepackage{tikz}
\usepackage{pgfplots}
\pgfplotsset{compat=1.3}
\usepackage{enumitem}
\usepackage[hidelinks]{hyperref}
\usepackage{cleveref}
\crefname{subsection}{subsection}{subsections}
\numberwithin{equation}{section}

\newcommand{\D}{\mathrm{D}}

\newtheorem{theorem}{Theorem}[section]
\newtheorem{proposition}[theorem]{Proposition}
\newtheorem{lemma}[theorem]{Lemma}

\newtheorem{corollary}[theorem]{Corollary}
\theoremstyle{remark}

\newtheorem{algorithm}{Algorithm}

\begin{document}

\title[Minimal residual methods for fully nonlinear PDE]{Minimal residual discretization of a class of
       fully nonlinear elliptic PDE}

\author[D.~Gallistl, N.~T.~Tran]{Dietmar Gallistl \and Ngoc Tien Tran}
\thanks{The authors received funding from the European Union's Horizon 2020 research and innovation programme (project DAFNE grant agreement No.~891734, and project RandomMultiScales, grant agreement No.~865751).}
\address[D.~Gallistl]{Friedrich-Schiller-Universit\"at Jena, 07743 Jena, Germany}
\email{dietmar.gallistl [at] uni-jena.de}
\address[N.~T.~Tran]{Universit\"at Augsburg, 86159 Augsburg, Germany}
\email{ngoc1.tran [at] uni-a.de}



\begin{abstract}
This work introduces finite element methods for a class of elliptic
fully nonlinear partial differential equations.
They are based on a minimal residual principle that builds upon
the Alexandrov--Bakelman--Pucci estimate.
Under rather general structural assumptions on the operator,
convergence of $C^1$ conforming and discontinuous Galerkin methods
is proven in the $L^\infty$ norm.
Numerical experiments on the performance of adaptive mesh refinement
driven by local information of the residual in two and three 
space dimensions are provided.
\end{abstract}

\maketitle

\section{Introduction}

Given an open bounded polyhedral Lipschitz domain $\Omega \subset \mathbb{R}^n$ and a right-hand side $f \in L^n(\Omega)$, we are interested in the numerical approximation
of strong solutions $u \in C(\overline{\Omega}) \cap W^{2,n}_\mathrm{loc}(\Omega)$
to second-order partial differential equations (PDE) in nondivergence form
\begin{align}\label{def:PDE}
	F(x, u, \nabla u, \D^2 u) = f \text{ a.e.~in } \Omega \quad\text{and}\quad u = g \text{ on } \partial \Omega
\end{align}
where the elliptic operator $F$ is \emph{proper} in the sense of structural assumption
\eqref{ineq:structure-assumption} below.
For a continuous operator $F$ and right-hand side $f$, strong solutions can be 
understood in the viscosity sense,
and finite difference methods (FDM) can preserve the maximum principle on the
discrete level and lead to a monotone scheme
with well-understood convergence theory \cite{BarlesSouganidis1991}.
In the same context,
also monotone finite element methods (FEM) \cite{JensenSmears2013,Jensen2017}
could be designed for certain classes of Hamilton--Jacobi--Bellman (HJB) equations.
On the other hand, monotone methods
have certain limitations in several regards and 
nonmonotone methods might offer better flexibility with respect to the approximation
order and mesh refinement. However, the design of FEMs for 
\eqref{def:PDE} turns out rather challenging.
The results in the literature mainly concern HJB equations
\begin{align*}
	0 = F(x,\D^2 u) \coloneqq \sup_{\alpha \in \mathcal{A}} (-A^\alpha : \D^2 u + \xi^\alpha) 
   \quad\text{in } \Omega.
\end{align*}
If the domain $\Omega$ is convex and $A^\alpha$ satisfies the so-called Cordes condition
(plus continuity assumptions),
then there exists a unique strong solution $u \in H^2(\Omega)$ to \eqref{def:PDE} for any $f \in L^2(\Omega)$ and $g \in H^{3/2}(\Omega)$.
This condition can be understood in the sense that the linearization
of \eqref{def:PDE} is sufficiently close to the Laplace equation.
The variational setting proposed in \cite{SmearsSueli2013,SmearsSueli2014} allows for access to finite element approximation \cite{SmearsSueli2014,GallistlSueli2019} with convergent adaptive mesh-refining algorithm \cite{GallistlSueli2019,KaweckiSmears2022}.
In two space dimensions, any uniformly elliptic operator 
without lower-order terms 
satisfies the Cordes condition, but the restriction becomes fairly strict for $n \geq 3$.

The goal of this paper is the design of convergent minimal residual methods 
for PDE in nondivergence form in extension to \cite{Tran2024}.
In contrast to the finite element schemes outlined above, solely local regularity of strong solutions is required for the convergence of the proposed FEM, which allows for nonsmooth domains with corner singularities.
In what follows we will write 
$$
  F[v] \coloneqq F(\bullet, v(\bullet), \nabla v(\bullet), \D^2 v(\bullet))
$$
and assume that there exists a strong solution $u \in W^{2,n}_\mathrm{loc}(\Omega)$.
The point of departure is the 
Alexandrov--Bakelman--Pucci (ABP) maximum principle for
the nonlinear Pucci operator, which implies the stability estimate
\begin{align}\label{ineq:stability}
	\|v - u\|_{L^\infty(\Omega)} \leq \Psi(v)
	\quad\text{for any }
	v \in W^{2,n}(\Omega)
\end{align}
with the residual
\begin{align}\label{def:Psi}
	\Psi(v)\coloneqq 
	\|g - v\|_{L^\infty(\partial \Omega)} + C\|f - F[v]\|_{L^n(\Omega)}
\end{align}
where $C$ is a generic constant.
We note that
$F[v] \in L^n(\Omega)$ holds under standard assumptions \eqref{ineq:structure-assumption}--\eqref{assumption:scaling} on $F$ below.
The right-hand side of \eqref{ineq:stability}
only depends on given data
and solving \eqref{def:PDE} becomes minimizing
$\Psi$, which is a suitable setting for stable 
subspace discretizations. The boundary conditions require
particular care.
For instance, let $V_h \subset W^{2,n}(\Omega)$ denote a finite dimensional 
subspace of $W^{2,n}(\Omega)$.
If $g$ is the trace of a discrete function, then the minimization of $\Psi$ among $v_h \in \{w_h \in V_h : w_h = g_h \text{ on } \partial \Omega\}$ is possible.
An example from \cite{Tran2024} shows that this finite element scheme may fail to approximate the exact solution $u$ unless $u \in W^{2,n}(\Omega)$ is sufficiently smooth.
However, this would require rather severe restrictions on the smoothness of the operator $F$ or the boundary $\partial \Omega$.
Instead, the minimization of $\Psi$ is carried out in the whole discrete space $V_h$.
We establish that, if
\begin{align}\label{ineq:minimization}
	\inf_{v \in W^{2,n}(\Omega)} \Psi(v) = 0,
\end{align}
then the sequence of discrete approximations converges uniformly to $u$.
This concept also leads to the design of convergent nonconforming FEM with well-understood averaging techniques.
A~posteriori error estimates are built-in thanks to \eqref{ineq:stability} and motivate an adaptive mesh-refining algorithm.
A (nonexhaustive) list of application satisfying \eqref{ineq:minimization} is presented and includes examples such as the Pucci or a regularized Monge--Amp\`ere equation.
Unless $F$ is linear, the residual $\Psi$ will be nonsmooth and nonconvex.

In comparison with the nonmonotone FEM of \cite{SmearsSueli2014}
for HJB equations,
our approach is built on a different solution concept.
The method of \cite{SmearsSueli2014} requires convexity of the domain
and the Cordes condition, but is superior to our method when these
conditions are met. It is based on the $W^{2,2}$ norm.
However, a general solution theory in $W^{2,2}$ in the absence 
of the Cordes condition or on non-convex domains is not known to us.
Our approach is based on strong solutions (which are at the same time
viscosity solutions) and requires only the structure conditions
around \eqref{ineq:structure-assumption}.
It provides error bounds in the $L^\infty$ norm, even is the 
solution is not bounded in the $W^{2,n}$ norm.

\medskip
The remaining parts of this paper are organized as follows. \Cref{sec:FEM} 
establishes stability estimates in the $L^\infty$ norm, which motivate the subsequent minimal residual methods and convergence of conforming FEM. \Cref{sec:NC-FEM} extends the results for conforming to nonconforming methods,
which are necessary for devising a practical method for $n\geq3$.
Applications covered by the theory are provided in
\Cref{sec:applications}.
Numerical experiments in \Cref{sec:numerical-examples} conclude this paper.

\medskip
Standard notation applies throughout this article.
The notation $\mathbb{S} \subset \mathbb{R}^{n \times n}$ denotes the subset of all symmetric matrices.
For $A,B \in \mathbb{S}$, the relation $A \geq 0$ means that $A$ is positive semidefinite and $A \leq B$ abbreviates $B - A \geq 0$.
The Pucci operators $\mathcal{P}^-_{\lambda,\Lambda}, \mathcal{P}^+_{\lambda,\Lambda} : \mathbb{S} \to \mathbb{R}$ associated with the ellipticity constants $0 < \lambda \leq \Lambda < \infty$ are defined, for any $M \in \mathbb{S}$, by
\begin{align}\label{def:pucci}
	\mathcal{P}_{\lambda,\Lambda}^-(M) \coloneqq \inf_{\lambda \mathrm{I} \leq A \leq \Lambda \mathrm{I}} (-A : M) \quad\text{and}\quad \mathcal{P}_{\lambda,\Lambda}^+(M) \coloneqq \sup_{\lambda \mathrm{I} \leq A \leq \Lambda \mathrm{I}} (-A : M).
\end{align}
Given any function $v : \Omega \to \mathbb{R}$, $v^+ \coloneqq \max\{0,v\}$ (resp.~$v^- \coloneqq \max\{0,-v\}$) denotes the positive (resp.~negative) part of $v$.

\section{Stability and conforming discretization}
\label{sec:FEM}

\subsection{Structure assumptions}
Throughout, we 
assume that there exist constants 
$0 < \lambda \leq \Lambda < \infty$, $\gamma \geq 0$, and $\mu \geq 0$ such that
\begin{align}\label{ineq:structure-assumption}
	\begin{split}
		\mathcal{P}^-_{\lambda,\Lambda}(M - N) - \gamma|p - q| - \mu(s-r)^+ \leq F(x,r,p,M) - F(x,s,q,N)&\\
		\leq \mathcal{P}^+_{\lambda,\Lambda}(M - N) + \gamma|p - q| + \mu(r-s)^+&
	\end{split}
\end{align}
holds for any $M, N \in \mathbb{S}$, $p,q \in \mathbb{R}^n$, $r,s \in \mathbb{R}$ with $|r|,|s| \leq R$, and a.e.~$x \in \Omega$.
This is a standard structure assumption on $F$ in the context of viscosity solutions \cite{CrandallIshiiLions1992,CaffarelliCrandallKocanSwiech1996} also known under the label of \emph{proper}, i.e., $F$ is uniformly elliptic and non-decreasing in $r$.
The dependence on $r$ in \eqref{ineq:structure-assumption} is simplified in comparison to the condition (SC) in \cite{CaffarelliCrandallKocanSwiech1996}, but can be extended to that case as well.
Furthermore, we assume that
\begin{align}\label{assumption:scaling}
	F(x,0,0,0) \in L^n(\Omega).
\end{align}
The combination of \eqref{ineq:structure-assumption}--\eqref{assumption:scaling} shows that $F$ is measurable in $(x,r,p,M)$ because $F$ is a Carath\'eodory function. In particular, 
$F[v]$ is a Lebesgue measurable function for any $v \in W^{2,1}_\mathrm{loc}(\Omega)$.

\subsection{Stability estimate}
This subsection establishes the stability estimate \eqref{ineq:stability}.
The fundamental tool is the ABP maximum principle
from \cite{CaffarelliCrandallKocanSwiech1996}.
\begin{lemma}[ABP maximum principle]\label{lem:maximum-principle}
	Let $f \in L^n(\Omega)$ and $\gamma > 0$. There exists a constant $C$ depending only on $\lambda$, $\Lambda$, $\gamma$, $n$, and $\mathrm{diam}(\Omega)$ such that (a)--(b) hold.
	\begin{enumerate}
		\item[(a)] If $v \in C(\overline{\Omega}) \cap W^{2,n}_\mathrm{loc}(\Omega)$ satisfies $f \leq \mathcal{P}^+_{\lambda,\Lambda}(\D^2 v) + \gamma|\nabla v| \quad\text{a.e.~in } \{v < 0\}$, then $\sup_\Omega v^- \leq \sup_{\partial \Omega} v^- + C\|f^-\|_{L^n(\Omega)}$.
		\item[(b)] If $v \in C(\overline{\Omega}) \cap W^{2,n}_\mathrm{loc}(\Omega)$ satisfies $\mathcal{P}^-_{\lambda,\Lambda}(\D^2 v) - \gamma|\nabla v| \leq f \quad\text{a.e.~in } \{0 < v\}$, then $\sup_\Omega v \leq \sup_{\partial \Omega} v + C\|f^+\|_{L^n(\Omega)}$.
	\end{enumerate}
\end{lemma}
\begin{proof}
	This result can be found in \cite[Proposition 3.3]{CaffarelliCrandallKocanSwiech1996}, where we used that $L^n$ strong supersolution (resp.~subsolution) are $L^n$ viscosity supersolution (resp.~subsolution) \cite[Lemma 2.5]{CaffarelliCrandallKocanSwiech1996}.
	Since the focus of this paper is not on the theory of $L^p$ viscosity solutions, we refer to \cite{CaffarelliCrandallKocanSwiech1996} for more details.
\end{proof}
\begin{theorem}[stability]\label{thm:stability}
	Suppose that $F$ satisfies \eqref{ineq:structure-assumption}--\eqref{assumption:scaling}.
	Then any $u, v \in C(\overline{\Omega}) \cap W^{2,n}_\mathrm{loc}(\Omega)$ with $F[u], F[v] \in L^n(\Omega)$ satisfy \eqref{ineq:stability}.
	The constant $C$ solely depends on $\lambda, \Lambda$, $\mu$, $n$, and $\mathrm{diam}(\Omega)$.
\end{theorem}
\begin{proof}
	The proof is essentially contained in \cite{CaffarelliCrandallKocanSwiech1996}. 
	The details are provided below for the sake of completeness.
	The structure assumption \eqref{ineq:structure-assumption} implies
	that, a.e.\ in $\Omega$, we have
	\begin{align}\label{ineq:proof-stability}
		\begin{split}
			&\mathcal{P}^-_{\lambda,\Lambda}(\D^2 (u - v)) - \gamma|\nabla(u - v)| - \mu(v - u)^+\\
			&\qquad\leq F[u] - F[v]
	  \leq \mathcal{P}^+_{\lambda,\Lambda}(\D^2 (u - v)) + \gamma|\nabla(u - v)| + \mu(u - v)^+.
		\end{split}
	\end{align}
	This proves
	$f - F[v] \leq \mathcal{P}^+_{\lambda,\Lambda}(\D^2 (u - v)) + \gamma|\nabla(u - v)|$ 
	a.e.~in $\{u-v < 0\}$ 
	and $\mathcal{P}^-_{\lambda,\Lambda}(\D^2 (u - v)) - \gamma|\nabla(u - v)| \leq f - F[v]$ a.e.~in $\{0 < u - v\}$,
    where we write $f=F[u]$.
	Hence, the ABP maximum principle for the Pucci operators from \Cref{lem:maximum-principle} concludes the proof.
\end{proof}

An immediate consequence of \Cref{thm:stability} is the uniqueness of strong solutions. Furthermore, if $u \in C(\overline{\Omega}) \cap W^{2,n}_\mathrm{loc}(\Omega)$ is a strong solution to \eqref{def:PDE}, then \Cref{thm:stability} provides \eqref{ineq:stability}.
We cannot expect that the functional $\Psi$
defined in \eqref{def:Psi} is efficient in the sense that $\Psi(v) \lesssim \|u - v\|_{L^\infty(\Omega)}$.
An upper bound is provided locally by error in the $W^{2,n}$ norm.
\begin{lemma}[local Lipschitz continuity]\label{lem:Lipschitz}
	Suppose that $F$ satisfies \eqref{ineq:structure-assumption}--\eqref{assumption:scaling}. Then any $v,w \in W^{2,n}(\omega)$ in an open subset $\omega \subset \Omega$ satisfy
	\begin{align*}
		\|F[v] - F[w]\|_{L^n(\omega)} \leq C_1\|v - w\|_{W^{2,n}(\omega)} .
	\end{align*}
	The constant $C_1$ solely depends on $\Lambda$, $\gamma$, $\mu$, $n$, and $\mathrm{diam}(\Omega)$.
\end{lemma}
\begin{proof}
	Since $-\sqrt{n}\Lambda|M|\leq \mathcal{P}^-_{\lambda,\Lambda}(M) \leq \mathcal{P}^+_{\lambda,\Lambda}(M) \leq \sqrt{n}\Lambda|M|$, \eqref{ineq:proof-stability} implies
	\begin{align*}
		|F[v] - F[w]|
		\leq \sqrt{n}\Lambda|\D^2(v - w)| + \gamma|\nabla(v - w)| + \mu|v - w|
	\end{align*}
	a.e.~in $\omega$. The left-hand side is a $L^n$ function and therefore, the claim follows.
\end{proof}
A first consequence of \Cref{lem:Lipschitz} is that $F[v] \in L^n(\omega)$ for any $v \in W^{2,n}(\omega)$ by choosing $w \coloneqq 0$ in \Cref{lem:Lipschitz}.
A second consequence of \Cref{lem:Lipschitz} is the continuity of $\Psi : W^{2,n}(\Omega) \to \mathbb{R}$. In fact, the inverse triangle inequality, \Cref{lem:Lipschitz}, and the Sobolev embedding $W^{2,n}(\Omega) \subset C(\overline{\Omega})$ provide, for any $v, w \in W^{2,n}(\Omega)$, that
\begin{align}\label{ineq:continuity}
	|\Psi(v) - \Psi(w)| \lesssim \|v - w\|_{L^\infty(\partial \Omega)} + \|v - w\|_{W^{2,n}(\Omega)} \lesssim \|v - w\|_{W^{2,n}(\Omega)}.
\end{align}

\subsection{Conforming discretization}
The convergence of conforming discretizations
relies on the assumption that \eqref{ineq:minimization} 
holds, which can be verified in the examples given in \Cref{sec:applications}.
Notice that the infimum cannot be attained as a minimum
unless the exact solution 
$u$ is smooth in the sense that it belongs to $W^{2,n}(\Omega)$.
Let $(V_j)_j$ denote a sequence of finite dimensional 
subspaces of $W^{2,n}(\Omega)$ with the density property 
$\lim_{j \to \infty} \min_{v_j \in V_j} \|v - v_j\|_{W^{2,n}(\Omega)} =0$ 
for all $v \in W^{2,n}(\Omega)$.
In other words, any $v \in W^{2,n}(\Omega)$ can be approximated by 
a sequence $(v_j)_j$ of discrete functions $v_j \in V_j$.
We think of the Argyris or Bogner--Fox--Schmid finite element space 
\cite{BrennerScott2008} on uniformly refined triangulations of the domain
$\Omega$.
Recall the residual $\Psi$ from \eqref{ineq:stability}--\eqref{def:Psi}.
The conforming method minimizes $\Psi$ in $V_j$.

\begin{proposition}[existence of discrete minimizers]\label{prop:existence-of-minimizers}
	The minimum of $\Psi$ in any finite dimensional 
	subspace $V_h \subset W^{2,n}(\Omega)$ is attained.
\end{proposition}
\begin{proof}
Since $\Psi$ is a nonnegative functional, we have $\inf\Psi(V_h)\geq0$.
	Let $(v_{\ell})_\ell \subset V_h$ be an infimizing sequence.
	In particular, $(\Psi(v_\ell))_\ell$ is uniformly bounded.
	\Cref{thm:stability} (with $u \coloneqq 0$) states
	\begin{align*}
		\|v_\ell\|_{L^\infty(\Omega)} 
		\leq \|v_\ell\|_{L^\infty(\partial \Omega)} 
		+ C\|F[v_\ell]\|_{L^n(\Omega)}.
	\end{align*}
	This and the triangle inequality provide 
	the bound 
	$$
	\|v_\ell\|_{L^\infty(\Omega)} \leq \|g\|_{L^\infty(\partial \Omega)} 
	+ C\|f\|_{L^n(\Omega)} + \Psi(v_\ell).
	$$
	Hence, $(v_\ell)_\ell$ is a bounded sequence in the finite dimensional space $V_h$
	and there exist $u_h \in V_h$ such that, up to some subsequence, 
	$\lim_{\ell \to \infty} \|u_h - v_\ell\|_{W^{2,n}(\Omega)} = 0$.
	This and the continuity of $\Psi$ in $V_h$ from \eqref{ineq:continuity}
	conclude the proof.
\end{proof}

The convergence of FEM is a consequence of the density of $\cup_j V_j$ in $W^{2,n}(\Omega)$.

\begin{theorem}[convergence of FEM]\label{thm:sufficiency}
	Suppose that $F$ satisfies \eqref{ineq:structure-assumption}--\eqref{assumption:scaling}. 
	If \eqref{ineq:minimization} holds,
	then $\lim_{j \to \infty} \min \Psi(V_j) = 0$.
	If, additionally, there exists a strong solution $u \in C(\overline{\Omega}) \cap W^{2,n}_\mathrm{loc}(\Omega)$ to \eqref{def:PDE},
	then the sequence of discrete minimizers $u_j \in \arg\min_{v_j \in V_j} \Psi(v_j)$ converges uniformly to $u$.
\end{theorem}
\begin{proof}
	Given $\varepsilon > 0$, the assumption \eqref{ineq:minimization} implies that
	that $\Psi(v) \leq \varepsilon$ for some $v \in W^{2,n}(\Omega)$.
	Let $v_j \in V_j$ denote the best-approximation of $v$ 
	in $V_j$ with respect to the $W^{2,n}$ norm.
	By the approximation property of $V_j$, $\lim_{j \to \infty} \|v - v_j\|_{W^{2,n}(\Omega)} = 0$.
	This and the continuity of $\Psi$ from \eqref{ineq:continuity} prove $\lim_{j \to \infty} \Psi(v_j) = \Psi(v) \leq \varepsilon$.
	Since $\Psi(u_j) \leq \Psi(v_j)$ for any $j$, we obtain $0 \leq \liminf_{j \to \infty} \Psi(u_j) \leq \limsup_{j \to \infty} \Psi(u_j) \leq \varepsilon$ for any $\varepsilon > 0$.
	This implies $\lim_{j \to \infty} \Psi(u_j) = 0$.
	If a strong solution $u$ exists, then the uniform convergence of $(u_j)_j$ to $u$ follows from \Cref{thm:stability}.
\end{proof}
Under additional smoothness assumptions on $u$, convergence rates can be derived from the approximation quality of the discrete space $V_j$.
\begin{corollary}[a~priori]\label{cor:a-priori}
	Suppose that $F$ satisfies \eqref{ineq:structure-assumption}--\eqref{assumption:scaling}. Let $u \in W^{2,n}(\Omega)$ be a strong solution to \eqref{def:PDE}.
	Then the sequence of discrete approximations $u_j \in \arg\min_{v_j \in V_j} \Psi(v_j)$ satisfies
	\begin{align*}
		\|u - u_j\|_{L^\infty(\Omega)} \lesssim \min_{v_j \in V_j} \|u - v_j\|_{W^{2,n}(\Omega)}.
	\end{align*}
\end{corollary}
\begin{proof}
	From \eqref{ineq:continuity} and $\Psi(u) = 0$, we deduce that $\Psi(v_j) \lesssim \|v - v_j\|_{W^{2,n}(\Omega)}$ for any $v_j \in V_j$. 
	This and \Cref{thm:stability} conclude the proof.
\end{proof}
We note that the a~priori result from \Cref{cor:a-priori} is only of limited interest because, if $u \in W^{2,n}(\Omega)$ is known a~priori, the minimal residual method can be organized more efficiently by prescribing the boundary data of discrete functions.
The regularity $u \in W^{2,n}(\Omega)$ implies that $g \in W^{2,n}(\Omega)$.
Define $W \coloneqq \{v \in W^{2,n}(\Omega) : v = 0 \text{ on } \partial \Omega\}$, $W_j \coloneqq \{v_j \in V_j: v_j = 0 \text{ on } \partial \Omega\}$, and
\begin{align*}
	\Phi(v) \coloneqq \|f - F[v]\|_{L^n(\Omega)} \quad\text{for any } v \in W^{2,n}(\Omega).
\end{align*}
\begin{proposition}[approximation for smooth solution]
	Suppose that $F$ satisfies \eqref{ineq:structure-assumption}--\eqref{assumption:scaling} and that a strong solution $u \in W^{2,n}(\Omega)$ to \eqref{def:PDE} exists.
	Furthermore, assume that $\lim_{j \to \infty} \min_{v_j \in W_j} \|v - v_j\|_{W^{2,n}(\Omega)} = 0$ for any $v \in W$.
	Let $g_j \in V_j$ approximate $g$ with $\lim_{j \to \infty} \|g - g_j\|_{W^{2,n}(\Omega)} = 0$.
	Then the sequence of discrete minimizers $u_j \coloneqq \arg\min \Phi(g_j + W_j)$ converges uniformly to $u$ as $j \to \infty$.
\end{proposition}
\begin{proof}
	Let $w_j$ denote the best approximation of $u - g$ onto $W_j$. Since $\lim_{j \to \infty} \|g - g_j\|_{W^{2,n}(\Omega)} = 0$ and $\lim_{j \to \infty} \min_{v_j \in W_j} \|v - v_j\|_{W^{2,n}(\Omega)} = 0$, $\lim_{j \to \infty} \|u - v_j\|_{W^{2,n}(\Omega)} = 0$ with $v_j \coloneqq w_j + g_j$. Hence, the continuity of $\Phi$ in $W^{2,n}(\Omega)$ from \eqref{ineq:continuity} leads to $\lim_{j \to \infty} \Phi(v_j) = \Phi(u) = 0$. This, $\Phi(u_j) \leq \Phi(v_j)$ for any $j$, and \Cref{thm:stability} concludes the proof.
\end{proof}

\subsection{Practical conforming FEM}\label{sub:cFEM}
In practice, the minimization of $\Psi$ from \eqref{def:Psi} in $V_j$ is reformulated as a constrained minimization problem to handle the nonsmooth part $\|g - v\|_{L^\infty(\partial \Omega)}$ \cite{Tran2024}. The relevant details are highlighted below. Let $g_j \in V_j$ be an approximation of $g$ on $\partial \Omega$. We assume that the $L^\infty(\partial \Omega)$ norm can be discretized in $V_j$, i.e., there exists a finite set $\mathcal{L}_j^b \subset \partial \Omega$ of points on the boundary with
\begin{align}\label{ineq:Linf-boundary}
	\|v_j\|_{L^\infty(\partial \Omega)} \leq c_\mathrm{eq}\max_{z \in \mathcal{L}_j^b} |v_j(z)|
\end{align} 
for any $v_j \in V_j$ with a constant $c_\mathrm{eq}$ independent of the index $j$. This assumption arises from the equivalence of norms in discrete space and holds if the traces of functions in $V_j$ are piecewise polynomials. For example, if $V_j$ is the BFS finite element space in two space dimensions, then $\mathcal{L}_j^b$ is the set of all Lagrange points associated with piecewise $P_3$ functions on the boundary.
Define the set
\begin{align*}
	\mathcal{A}_j(g_j) \coloneqq \{(t,v_j) \in \mathbb{R}_{\geq 0} \times V_j : - t \leq g_j(z) - v_j(z) \leq t \text{ for any } z \in \mathcal{L}_j^b\}
\end{align*}
of discrete admissible functions.
Introducing an additional variable $t$ for the non\-smooth part,
the practical method seeks a minimizer $(s_j,u_j) \in \mathcal{A}_j(g_j)$ of
\begin{align}\label{def:Phi-pract}
	\Phi_j(t,v_j) \coloneqq \sigma t^n + \|f - F[v_j]\|^n_{L^n(\Omega)} \quad\text{among } (t,v_j) \in \mathcal{A}_j(g_j)
\end{align}
with a fixed positive parameter $\sigma > 0$.
It is proven in \cite{Tran2024} that the minimization of $\Psi$ in $V_j$ and of $\Phi_j$ in $\mathcal{A}_j(g_j)$ are equivalent in the sense that
\begin{align*}
	\Phi_j(s_j,u_j) \approx \min_{v_j \in V_j} \big(\|g_j - v_j\|_{L^\infty(\partial \Omega)}^n + \|f - F[v_j]\|_{L^n(\Omega)}^n\big)
\end{align*}
with $s_j = \max_{z \in \mathcal{L}_j^b} |g_j(z) - u_j(z)|$.
The constants hidden in the notation $\approx$ are independent of $j$. This leads to the following convergence result.
\begin{theorem}[convergence of practical FEM]
	Suppose that $F$ satisfies \eqref{ineq:structure-assumption}--\eqref{assumption:scaling} and $\lim_{j \to \infty} \|g - g_j\|_{L^\infty(\partial \Omega)} = 0$.
	If \eqref{ineq:minimization} and \eqref{ineq:Linf-boundary} hold,
	then
	\begin{align*}
		\lim_{j \to \infty} \min \Phi_j(\mathcal{A}_j(g_j)) = 0
	\end{align*}
	If, additionally, there exists a strong solution $u \in C(\overline{\Omega}) \cap W^{2,n}_\mathrm{loc}(\Omega)$ to \eqref{def:PDE},
	then the sequence of discrete minimizers $(s_j,u_j) \in \arg\min \Phi_j(\mathcal{A}_j(g_j))$ satisfies $\lim_{j \to \infty} \|u - u_j\|_{L^\infty(\Omega)} = 0$.
\end{theorem}
\begin{proof}
	Since the arguments from the linear case \cite{Tran2024} carry over, further details are omitted.
\end{proof}

We note that the two additive terms in the definition of $\Phi_j$
in \eqref{def:Phi-pract} can in principle be balanced with
appropriate weights. We disregard this possibility in our qualitative
convergence analysis.

\section{Nonconforming FEM}\label{sec:NC-FEM}
In more than two space dimensions,
conforming methods are difficult to implement due to the large number of local degrees of freedom.
This section proposes a nonconforming FEM with the piecewise polynomial trial 
space $V_\mathrm{nc}(\mathcal{T}_j) \coloneqq P_k(\mathcal{T}_j)$,
that is, piecewise polynomial functions of degree at most $k \geq 2$ on 
a uniformly shape-regular sequence $(\mathcal{T}_j)_j$ of simplicial triangulations 
of $\Omega$ such that the maximal mesh-size 
$h_j \coloneqq \max_{T \in \mathcal{T}_j} \mathrm{diam}(T)$ vanishes in the limit as $j \to \infty$.
For piecewise smooth functions, $\nabla_\mathrm{pw}$ and 
$\D^2_\mathrm{pw}$ denote the piecewise gradient and Hessian,
respectively, without explicit reference to the triangulation $\mathcal{T}_j$.
Let $\mathcal{F}_j(\Omega)$ denote the set of interior faces
(the $(n-1)$-dimensional hyperfaces) of $\mathcal{T}_j$.
For any $p \in (1, \infty)$, define the stabilization $s_j : V_\mathrm{nc}(\mathcal{T}_j) \to \mathbb{R}$ by
\begin{align}\label{def:stabilization}
	s_j(v_j; p) \coloneqq \sum_{F \in \mathcal{F}_j(\Omega)} 
	 \left(h_F^{1-2p} \|[v_j]_F\|^p_{L^p(F)}
	 +
	 h_F^{1-p} \|[\nabla_\mathrm{pw} v_j]_F\|^p_{L^p(F)}\right).
\end{align}
Here, the brackets $[\cdot]_F$ denote the jump of a function across 
the face $F$.
We abbreviate $s_j(\bullet; n) \coloneqq s_j(v_j)$ and denote a piecewise
version of $F$ by
$$
 F_\mathrm{pw}[v] \coloneqq 
 F(\bullet, v(\bullet), \nabla_\mathrm{pw} v(\bullet), \D_\mathrm{pw}^2 v(\bullet)).
$$
The residual in the nonconforming case is defined as
\begin{align}\label{def:psi-nc}
	\begin{split}
		\Psi_j^{\mathrm{nc}}(v_j) &\coloneqq \|g - v_j\|_{L^\infty(\partial \Omega)}
		   + \|f - F_\mathrm{pw}[v_j]\|_{L^n(\Omega)} + s_j(v_j)^{1/n}.
	\end{split}
\end{align}

The analysis of nonconforming FEM relies on the construction of an enrichment operator.

\begin{lemma}[enrichment operator]\label{lem:enrichment}
    Let $n\in\{2,3\}$.
	Suppose that $k \geq 2$, then
	there exists a linear operator $\mathcal{J}_j : V_\mathrm{nc}(\mathcal{T}_j) \to W^{2,\infty}(\Omega)$ such that
	\begin{align}
			h_j^{-2p}\|v_j - \mathcal{J}_j v_j\|_{L^p(\Omega)}^p + h_j^{-p}\|\nabla_\mathrm{pw} (v_j - \mathcal{J}_j v_j)\|_{L^p(\Omega)}^p
			+ \|\D_\mathrm{pw}^2(v_j - \mathcal{J}_j v_j)\|_{L^p(\Omega)}^p
			\leq C_2 s_j(v_j;p)
		\label{ineq:local-enrichment-operator}
	\end{align}
	for all $v_j \in V_\mathrm{nc}(\mathcal{T}_j)$ and $p \in (1,\infty)$
	with a constant $C_2$ depending on $k$, $n$, $p$, and the shape regularity of $\mathcal{T}_j$.
\end{lemma}

\begin{proof}
	Averaging techniques with a local version of \eqref{ineq:local-enrichment-operator} are well understood in the literature \cite{BrennerGudiSung2010,GeorgoulisHoustonVirtanen2011,Gallistl2015}.
	The operator $\mathcal{J}_j$ maps into 
	a finite element space of continuously differentiable piecewise polynomial functions of degree at most $m \geq k$.
	The Hsieh–Clough–Tocher (HCT) or Worsey--Farin
	macro elements \cite{CloughTocher1965,WorseyFarin1987} are the target space of 
	choice because the degrees of freedom only depend on the values 
	of the function and its first derivative.
	They are (currently) available for arbitrary polynomial degree 
	$m \geq 3$ in two and three space dimensions 
	\cite{DouglasDupontPercellScott1979,GuzmanLischkeNeilan2022}.
	In our case, we use a generalization of operators proposed in
	\cite{NeilanWu2019,KaweckiSmears2021}.
	The precise assignment of degrees of freedom is described in
	\cite[Lemma~3.5]{GallistlTian2024}, with the difference that
	we do not impose the zero boundary conditions in the 
	definition of the operator $\mathcal{J}_j$.
	This is the reason why no boundary jumps occur on the right-hand side
	of \eqref{ineq:local-enrichment-operator}.
	Therefore a dependence on the interior angles of the polygonal
	domain as described in \cite{KaweckiSmears2021} is avoided.
\end{proof}

In the following, results from the conforming case are extended to the nonconforming case.
We start with an extension of the stability result from \Cref{thm:stability} to handle discontinuous functions,
as a motivation of the residual of \eqref{def:psi-nc}.

\begin{theorem}[stability for NC-FEM]\label{thm:stability-nc}
	Suppose that $F$ satisfies \eqref{ineq:structure-assumption}.
	Any $v \in C(\overline{\Omega}) \cap W^{2,n}_\mathrm{loc}(\Omega)$ with $F(\bullet, v, \nabla v, \D^2 v) \in L^n(\Omega)$ and $v_j \in V_\mathrm{nc}(\mathcal{T}_j)$ satisfies
	\begin{align*}
		\|v - v_j\|_{L^\infty(\Omega)} \leq \|v - v_j\|_{L^\infty(\partial \Omega)}
		+ C_3(\|F[v] - F_\mathrm{pw}[v_j]\|_{L^n(\Omega)}^n + s_j(v_j))^{1/n}.
	\end{align*}
	The constant $C_3$ solely depends on $\lambda$, $\Lambda$, $\mu$, $k$, $n$, $\mathrm{diam}(\Omega)$, and the shape regularity of $\mathcal{T}_j$.
\end{theorem}
\begin{proof}
	The triangle inequality reads as
	\begin{align}\label{ineq:proof-stability-nc-split}
		\|v - v_j\|_{L^\infty(\Omega)} \leq \|v - \mathcal{J}_j v_j\|_{L^\infty(\Omega)} + \|\mathcal{J}_j v_j - v_j\|_{L^\infty(\Omega)}.
	\end{align}
	Since $\mathcal{J}_j v_j \in W^{2,n}(\Omega)$, \Cref{thm:stability} and the triangle inequality provide
	\begin{align}\label{ineq:stability-triangle}
		\begin{split}
			\|v - \mathcal{J}_j v_j\|_{L^\infty(\Omega)} 
			\leq
			&\|v - v_j\|_{L^\infty(\partial \Omega)} + \|v_j - \mathcal{J}_j v_j\|_{L^\infty(\partial \Omega)}\\
			&\quad+ C\|F[v] - F_\mathrm{pw}[v_j]\|_{L^n(\Omega)}
			+ C\|F_\mathrm{pw}[v_j] - F[\mathcal{J}_j v_j]\|_{L^n(\Omega)}.
		\end{split}
	\end{align}
	The piecewise application of \Cref{lem:Lipschitz} and the bound \eqref{ineq:local-enrichment-operator} lead to
	\begin{align}\label{ineq:proof-stability-nc-L-continuity}
		\begin{split}
			\|F_\mathrm{pw}[v_j] - F[\mathcal{J}_j v_j]\|_{L^n(\Omega)}
			\lesssim \|v_j - \mathcal{J}_j v_j\|_{W^{2,n}(\mathcal{T}_j)} \lesssim s_j(v_j)^{1/n}.
		\end{split}
	\end{align}
	Furthermore, $\varphi \coloneqq v_j - \mathcal{J}_j v_j \in W^{2,n}(\mathrm{int}(T))$ for each cell $T \in \mathcal{T}_j$ and so, $\varphi \circ L \in W^{2,n}(T_\mathrm{ref})$ with the affine transformation 
	$L : T_\mathrm{ref} \to T$
	to some fixed reference domain $T_\mathrm{ref}$.
	The Sobolev embedding implies $\|\varphi\|_{L^\infty(T)} = \|\varphi \circ L\|_{L^\infty(T_\mathrm{ref})} \lesssim \|\varphi \circ L\|_{W^{2,n}(T_\mathrm{ref})}$.
	A scaling argument establishes
	$$
	\|\varphi \circ L\|_{W^{2,n}(T_\mathrm{ref})} \lesssim h_T^{-1}\|\varphi\|_{L^n(T)} + h_T^{n-1}\|\nabla \varphi\|_{L^n(T)} + h_T^{2n-1}\|\D^2 \varphi\|_{L^n(T)}.
	$$
	Hence, \Cref{lem:enrichment} provides
	$\|v_j - \mathcal{J}_j v_j\|_{L^\infty(\Omega)} \lesssim h_j s_j(v_j)$.
	The combination of this with \eqref{ineq:proof-stability-nc-split}--\eqref{ineq:proof-stability-nc-L-continuity} concludes the proof.
\end{proof}

The following results are analogous to those in the conforming case.
\begin{proposition}[existence of discrete minimizers NC-FEM]
	The minimum of $\Psi_j^\mathrm{nc}$ in $V_\mathrm{nc}(\mathcal{T}_j)$ is attained.
\end{proposition}
\begin{proof}
	The claim follows from similar arguments as those presented in the proof of \Cref{prop:existence-of-minimizers}.
\end{proof}

\begin{theorem}[convergence of NC-FEM]\label{thm:convergence-NC}
	Suppose that $F$ satisfies \eqref{ineq:structure-assumption}--\eqref{assumption:scaling}.
	If \eqref{ineq:minimization} holds,
	then $\lim_{j \to \infty} \min \Psi_j^\mathrm{nc}(V_\mathrm{nc}(\mathcal{T}_j)) = 0$.
	If, additionally, there exists a strong solution $u \in C(\overline{\Omega}) \cap W^{2,n}_\mathrm{loc}(\Omega)$ to \eqref{def:PDE},
	then the sequence of discrete minimizers $u_j \in \arg\min \Psi_j^\mathrm{nc}(V_\mathrm{nc}(\mathcal{T}_j))$ converges to $u$ with $\lim_{j \to \infty} \|u - u_j\|_{L^\infty(\Omega)} = 0$.
\end{theorem}

\begin{proof}
	Given $\varepsilon > 0$, the density of $C^{\infty}(\overline{\Omega})$ in $W^{2,n}(\Omega)$
	and assumption \eqref{ineq:minimization} imply that there exists 
	$v \in C^\infty(\overline{\Omega})$ such that $\Psi(v) \leq \varepsilon$.
	Let $v_j \coloneqq \Pi_{\mathcal{T}_j}^k v$ denote the $L^2$ projection of $v$ onto $V_\mathrm{nc}(\mathcal{T}_j)$.
	The approximation property of the $L^2$ projection 
	$\Pi_{\mathcal{T}_j}^k$
	implies that
	$\lim_{j \to \infty} \|v - v_j\|_{L^{\infty}(\Omega)} \lesssim h_j\|v\|_{W^{1,\infty}(\Omega)} \to 0$ 
	and $\|v - v_j\|_{W^{2,n}(\mathcal{T}_j)} \lesssim h_j\|v\|_{W^{3,n}(\Omega)} \to 0$ 
	as $j \to \infty$.
	This and a piecewise application of \Cref{lem:Lipschitz} lead to
	\begin{align}\label{eq:proof-convergence-nc-FEM-psi}
		\lim_{j \to \infty} (\|g - v_j\|_{L^\infty(\partial \Omega)} + \|f - F_\mathrm{pw}[v_j]\|_{L^n(\Omega)}) = \Psi(v).
	\end{align}
	The jumps $[\D^\alpha v]_F$ along any interior side $F \in \mathcal{F}_j(\Omega)$ and multiindex $\alpha$ vanish.
	Therefore, the trace inequality implies
	\begin{align}\label{ineq:stabilization-bound}
	s_j(v_j) \lesssim 
	h_j^{-2n}\|v - v_j\|_{L^n(\Omega)}^n
	+h_j^{-n}\|\nabla_\mathrm{pw}(v - v_j)\|_{L^n(\Omega)}^n
	+\|\D^2_\mathrm{pw}(v - v_j)\|_{L^n(\Omega)}^n.
	\end{align}
	By standard approximation properties of $v_j$ \cite[Lemma 11.18]{ErnGuermond2021},
	this vanishes in the limit as $j \to \infty$, 
	and \eqref{eq:proof-convergence-nc-FEM-psi} implies 
	$$
	0 \leq \limsup_{j \to \infty} \Psi_j^\mathrm{nc}(u_j) 
	  \leq \lim_{j \to \infty} \Psi_j^\mathrm{nc}(v_j) 
	   \leq \varepsilon.
	$$
	The parameter $\varepsilon > 0$ was chosen arbitrary, whence $\lim_{j \to \infty} \Psi_j^\mathrm{nc}(u_j) = 0$. \Cref{thm:stability-nc} concludes $\lim_{j \to \infty} \|u - u_j\|_{L^\infty(\Omega)} = 0$ if strong solutions exist.
\end{proof}

Under additional smoothness assumptions, convergence rates can be derived.

\begin{corollary}[a~priori NC-FEM]\label{cor:a-priori-nc}
    Let $n\in\{2,3\}$ and $k\geq2$.
	Suppose that $F$ satisfies \eqref{ineq:structure-assumption}--\eqref{assumption:scaling}. 
	Let $u \in W^{2,n}(\Omega)$ be a strong solution to \eqref{def:PDE}.
	Then the sequence of discrete approximations $u_j \in \arg\min \Psi(V_\mathrm{nc}(\mathcal{T}_j))$ satisfies
	$\|u - u_j\|_{L^\infty(\Omega)} \lesssim \|\D^2_\mathrm{pw}(1 - \Pi_{\mathcal{T}_j}^k) u\|_{L^n(\Omega)}$.
	In particular, if $u \in W^{2,n}(\Omega) \cap W^{k+1,n}(\mathcal{T}_j)$, then $\|u - u_j\|_{L^\infty(\Omega)} \lesssim h_j^{k-1}$.
\end{corollary}
\begin{proof}
	Let $v_j \coloneqq \Pi_{\mathcal{T}_j}^k u$ denote the $L^2$ projection of $u$ onto $P_k(\mathcal{T}_j)$.
	The piecewise application of the Sobolev embedding as outlined in the proof of \Cref{thm:stability-nc} leads to 
	\begin{align*}
	\|u - v_j\|_{L^\infty(\Omega)}
	&\lesssim h_j^{-1}\|u - v_j\|_{L^n(\Omega)} + h_j^{n-1}\|\nabla_\mathrm{pw}(u - v_j)\|_{L^n(\Omega)}
	+ h_j^{2n-1}\|\D^2(u - v_j)\|_{L^n(\Omega)} 
	\\
	&\lesssim h_j\|\D^2_\mathrm{pw}(u - v_j)\|_{L^n(\Omega)}.
	\end{align*}
	This, a piecewise application of \Cref{lem:Lipschitz}, and \eqref{ineq:stabilization-bound} result in $\Psi_j^\mathrm{nc}(v_j) \lesssim \|\D^2_\mathrm{pw}(u - v_j)\|_{L^n(\Omega)}$. The stability result from \Cref{thm:stability-nc} and $\|u - u_j\|_{L^\infty(\Omega)} \lesssim \Psi_j^\mathrm{nc}(u_j) \leq \Psi_j^\mathrm{nc}(v_j)$ conclude the proof.
\end{proof}

In analogy to the discussion of \Cref{sub:cFEM},
in practice, the minimization of $\Psi_j^\mathrm{nc}$ from \eqref{def:psi-nc} is reformulated as a constrained minimization problem to handle the nonsmooth part $\|g - v_j\|_{L^\infty(\partial \Omega)}$ \cite{Tran2024}. The relevant details are highlighted below. Let $g_j \in V_\mathrm{nc}(\mathcal{T}_j)$ be an approximation of $g$ on $\partial \Omega$. Then the equivalence of norms in finite dimensional spaces shows that
\begin{align*}
	c_\mathrm{eq}^{-1} \|g_j - v_j\|_{L^\infty(\Omega)} \leq \max_{T \in \mathcal{T}_j} \max_{z \in \mathcal{L}_j^b \cap T} |g_j(z) - v_j|_T(z)| \leq \|g_j - v_j\|_{L^\infty(\Omega)},
\end{align*}
where $\mathcal{L}_j^b$ denotes the set of Lagrange nodes associated with the discrete space $P_k(\mathcal{T}_j) \cap H^1(\Omega)$ on the boundary $\partial \Omega$. 
The constant $c_\mathrm{eq}$ depends on $k$ and the shape regularity but not on the mesh-size of $\mathcal{T}_j$.
Define the set
\begin{align*}
	\mathcal{A}_\mathrm{nc}(g_j, \mathcal{T}_j) \coloneqq \{(t,v_j) \in \mathbb{R}_{\geq 0} \times V_\mathrm{nc}(\mathcal{T}_j) : - t \leq g_j(z) - v_j(z) \leq t \text{ for any } z \in \mathcal{L}_j^b\}
\end{align*}
of discrete admissible functions.
Introducing an additional variable $t$ for the non\-smooth part,
the practical method seeks a minimizer $(s_j,u_j) \in \mathcal{A}_\mathrm{nc}(g_j,\mathcal{T}_j)$ of
\begin{align}\label{def:Phi-nc-pract}
	\Phi_j^\mathrm{nc}(t,v_j) \coloneqq \sigma t^n + \|f - F_\mathrm{pw}[v_j]\|^n_{L^n(\Omega)} + \tau s_j(v_j) \quad\text{among } (t,v_j) \in \mathcal{A}_\mathrm{nc}(g_j,\mathcal{T}_j)
\end{align}
with fixed positive parameters $\sigma, \tau > 0$.
It is proven in \cite{Tran2024} that the minimization of $\Psi^\mathrm{nc}_j$ in $V_\mathrm{nc}(\mathcal{T}_j)$ and of $\Phi_j^\mathrm{nc}$ in $\mathcal{A}_\mathrm{nc}(g_j,\mathcal{T}_j)$ are equivalent in the sense that
\begin{align*}
	\Phi_j^\mathrm{nc}(s,u_j) \approx \min_{v_j \in V_j} \big(\|g_j - v_j\|_{L^\infty(\partial \Omega)}^n + \|f + F_\mathrm{pw}[v_j]\|_{L^n(\Omega)}^n + s_j(v_j)\big).
\end{align*}
The constants hidden in the notation $\approx$ depend on $\sigma$, $\tau$, $k$, and the shape regularity of $\mathcal{T}_j$. This leads to the following convergence result.
\begin{theorem}[convergence of practical NC-FEM]
	Suppose that $F$ satisfies \eqref{ineq:structure-assumption}--\eqref{assumption:scaling} and $\lim_{j \to \infty} \|g - g_j\|_{L^\infty(\partial \Omega)} = 0$.
	If \eqref{ineq:minimization} holds,
	then
	\begin{align*}
		\lim_{j \to \infty} \min \Phi_j^\mathrm{nc}(\mathcal{A}_\mathrm{nc}(g_j,\mathcal{T}_j)) = 0.
	\end{align*}
	If, additionally, there exists a strong solution $u \in C(\overline{\Omega}) \cap W^{2,n}_\mathrm{loc}(\Omega)$ to \eqref{def:PDE},
	then the sequence of discrete minimizers $(s_j,u_j) \in \arg\min \Phi_j^\mathrm{nc}(\mathcal{A}_\mathrm{nc}(g_j,\mathcal{T}_j))$ satisfies $\lim_{j \to \infty} \|u - u_j\|_{L^\infty(\Omega)} = 0$.
\end{theorem}
\begin{proof}
	Since the arguments from the linear case \cite{Tran2024} carry over, further details are omitted.
\end{proof}

\section{Applications}\label{sec:applications}
This section provides a nonexhaustive overview of examples that
satisfy \eqref{ineq:minimization}.

\subsection{Linear PDE}
The linear case was already analyzed in \cite{Tran2024}. We recall the relevant result therein and refer to \cite{Tran2024} for further details.
\begin{proposition}
	Suppose that $A \in C(\overline{\Omega};\mathbb{S})$ with $\lambda \mathrm{I} \leq A \leq \Lambda \mathrm{I}$ pointwise in $\Omega$, $b \in L^\infty(\Omega; \mathbb{R}^n)$, and $0 \leq c \leq L^\infty(\Omega)$.
	Then for any $f \in L^n(\Omega)$ and $g \in C(\partial \Omega)$, there exists a unique strong solution $u \in C(\overline{\Omega}) \cap W^{2,n}_\mathrm{loc}(\Omega)$ 
	to the PDE \eqref{def:PDE} with 
	$F(x,r,p,M) \coloneqq - A(x):M + b(x) \cdot p + c(x)r$ and \eqref{ineq:minimization} holds.
\end{proposition}

\subsection{Pucci extremal operator}\label{sub:Pucci}
Recall $\mathcal{P}^+_{\lambda,\Lambda}(M) \coloneqq \sup_{\lambda \mathrm{I} \leq A \leq \Lambda \mathrm{I}} (-A : M)$. Given $\mu \geq 0$, define
\begin{align}\label{def:Pucci-equation}
	F(p,M) \coloneqq \mathcal{P}_{\lambda,\Lambda}^+(M) + \mu|p|
\end{align}
\begin{proposition}
	For any $f \in L^n(\Omega)$ and $g \in C(\partial \Omega)$, there exists a unique strong solution $u \in C(\overline{\Omega}) \cap W^{2,n}_\mathrm{loc}(\Omega)$ 
	to the PDE \eqref{def:PDE} with $F$ from \eqref{def:Pucci-equation} and \eqref{ineq:minimization} holds.
\end{proposition}
\begin{proof}
	The existence of strong solutions is established in \cite[Corollary 3.10]{CaffarelliCrandallKocanSwiech1996}. The proof of \eqref{ineq:minimization} follows the arguments given in the proof of \Cref{prop:HJB-strong-solution} below.
\end{proof}

\subsection{Hamilton--Jacobi--Bellman equations}
This class is of particular interest due to the existence of strong solutions arising from the Evans--Krylov theory \cite{Evans1982,Krylov1983}.

\subsubsection{Strong solutions}
For an index set $\mathcal{A}$, let $A^\alpha \in C(\overline{\Omega}; \mathbb{S})$ with $\lambda \mathrm{I}\leq A^\alpha \leq \Lambda \mathrm{I}$ pointwise in $\Omega$ for any $\alpha \in \mathcal{A}$ be given and set
\begin{align}\label{def:HJB-strong-solution}
	F(x,M) \coloneqq \sup_{\alpha \in \mathcal{A}} (-A^\alpha(x): M).
\end{align}

\begin{proposition}[HJB with strong solution]\label{prop:HJB-strong-solution}
	Suppose that $A^\alpha$, $\alpha \in \mathcal{A}$, is uniformly bounded in $C^{0,\beta}(\overline{\Omega})$ for a fixed $0 < \beta \leq 1$ independent of $\alpha$, i.e., there exists $C > 0$ such that $\|A^\alpha\|_{C^{0,\beta}(\overline{\Omega})} < C$.
	Then for any $f \in L^n(\Omega)$ and $g \in C(\partial \Omega)$, there exists a unique strong solution $u \in C(\overline{\Omega}) \cap W^{2,n}_\mathrm{loc}(\Omega)$ to the PDE \eqref{def:PDE} with $F$ from \eqref{def:HJB-strong-solution} and \eqref{ineq:minimization} holds.
\end{proposition}
The proof of \Cref{prop:HJB-strong-solution} requires the following result on H\"older continuity of nonlinear equations.
\begin{theorem}[global H\"older continuity]\label{thm:Hoelder}
	Suppose that $F$ satisfy \eqref{ineq:structure-assumption}--\eqref{assumption:scaling}.
	Given $f \in L^n(\Omega)$ and $g \in C^{0,\beta}(\partial \Omega)$,
	let $u \in C(\overline{\Omega}) \cap W^{2,n}_\mathrm{loc}(\Omega)$ be a strong solution to \eqref{def:PDE}.
	Then $u \in C^{0,\widetilde{\beta}}(\overline{\Omega})$ is H\"older continuous with an exponent $\widetilde{\beta}$ that solely depends on
	$n$, $\lambda$, $\Lambda$, $\beta$, and the exterior cone condition of $\Omega$.
	Furthermore,
	\begin{align*}
		|u(x) - u(y)| \leq C_4|x - y|^{\widetilde{\beta}} \quad\text{for any } x,y \in \overline{\Omega}.
	\end{align*}
	The constant $C_4$ exclusively depends on
	$n$, $\lambda$, $\Lambda$, $\gamma$, $\mu$, $\|f\|_{L^n(\Omega)}$, $\|g\|_{C^{0,\beta}(\partial \Omega)}$, $\mathrm{diam}(\Omega)$,
	and the cone condition of $\Omega$.
\end{theorem}
\begin{proof}
	A proof can be found, e.g., in \cite[Section 9.8--9.9]{GilbargTrudinger2001} or \cite[Theorem 6.2]{KoikeSwiech2009}.
\end{proof}
\begin{proof}[Proof of \Cref{prop:HJB-strong-solution}]
	The existence of solutions is proven in \cite[Example 3.11]{CaffarelliCrandallKocanSwiech1996}.
	The strategy from \cite{Tran2024} to verify \eqref{ineq:minimization}
	applies to all examples of this section and is only carried out in this case.
	For any $\delta > 0$, define the set $\Omega_\delta \coloneqq \{x \in \mathbb{R}^n : \mathrm{dist}(x,\overline{\Omega}) < \delta\}$. If $\delta$ is sufficiently small, then
	$\Omega_\delta$ is a Lipschitz domain so that the characteristics of the Lipschitz boundary of $\Omega_\delta$ and $\Omega$ coincide \cite{Doktor1976}.
	As a consequence, $\Omega_\delta$ satisfies a uniform exterior cone condition independent of $\delta$ for sufficiently small $\delta$ \cite{Grisvard2011}.
	Fix one of this $\delta$.
	The functions $A^\alpha$ are extended outside of $\Omega$ so that $A^\alpha \in C^{0,\beta}(\overline{\Omega}_\delta;\mathbb{S})$ and $\|A^\alpha\|_{C^{0,\beta}(\overline{\Omega}_\delta)} \leq C$ for a constant $C$ independent of $\alpha$ \cite[Corollary 1]{McShane1934}. We abuse the notation and use $A^\alpha$ for the extension as well.
	For $\varepsilon > 0$, let $\widetilde{g} \in C^\infty(\overline{\Omega}_\delta)$ be given such that $\|g - \widetilde{g}\|_{L^\infty(\partial \Omega)} \leq \varepsilon/2$.
	For any $j \in \mathbb{N}$, let $u_j \in C(\overline{\Omega}_{\delta/j}) \cap W^{2,n}_\mathrm{loc}(\Omega_{\delta/j})$ be the strong solution to
	\begin{align}\label{def:PDE-extension}
		F(x,\D^2 u_j) = f \text{ in } \Omega_{\delta/j} \quad\text{and}\quad u_j = \widetilde{g} \text{ on } \partial \Omega_{\delta/j}.
	\end{align}
	From \Cref{thm:Hoelder}, we deduce that there exists a constant $C_4$ independent of $j$ such that $|u_j(x) - u_j(y)| \leq C_4|x-y|^{\widetilde{\beta}}$ holds for any $j \in \mathbb{N}$ and $x,y \in \overline{\Omega}_{\delta/j}$.
	Thus, $|u_j(x) - u_j(x')| \leq C_4|x - x'|^{\widetilde{\beta}} = C_4(\delta/j)^{\widetilde{\beta}}$ for any $x \in \partial \Omega$ and its best-projection $x'$ onto $\partial \Omega_{\delta/j}$.
	Since $\widetilde{g}$ is Lipschitz continuous, this, a triangle inequality, and $u_j(x') = \widetilde{g}(x')$ prove
	\begin{align*}
		\|\widetilde{g} - u_j\|_{L^\infty(\partial \Omega)} \leq \max_{x \in \partial \Omega} (|\widetilde{g}(x) - \widetilde{g}(x')| + |u_j(x) - u_j(x')|) \leq L(\delta/j) + C_4(\delta/j)^{\widetilde{\beta}},
	\end{align*}
	where $L$ is the Lipschitz constant of $\widetilde{g}$.
	By choosing $j$ sufficiently large, we obtain $\|\widetilde{g} - u_j\|_{L^\infty(\partial \Omega)} \leq \varepsilon/2$.
	The combination of this with \eqref{def:PDE-extension} and a triangle inequality concludes $\Psi(u_j) = \|g - u_j\|_{L^\infty(\partial \Omega)} \leq \|g - \widetilde{g}\|_{L^\infty(\partial \Omega)} + \|\widetilde{g} - u_j\|_{L^\infty(\partial \Omega)} \leq \varepsilon$.
	Since $\varepsilon$ was chosen arbitrary and $u_j|_\Omega \in W^{2,n}(\Omega)$, this concludes the proof. 
\end{proof}

\subsubsection{Classical solutions}\label{sub:HJB-classical}
If the index set $\mathcal{A}$ is countable, we have classical solutions for H\"older continuous coefficients.
\begin{proposition}[HJB with classical solution]\label{prop:HJB-classical-solution}
	Suppose that $\mathcal{A}$ is a countable set.
	Let $A^\alpha \in C^{0,\beta}(\overline{\Omega};\mathbb{S})$, $b^\alpha \in C^{0,\beta}(\overline{\Omega};\mathbb{R}^n)$, $0 \leq c^\alpha \in C^{0,\beta}(\overline{\Omega})$, $\xi^\alpha \in C^{0,\beta}(\overline{\Omega})$ for any $\alpha \in \mathcal{A}$ and a fixed $0 < \beta \leq 1$ (independent of $\alpha$) be given such that
	their $C^{0,\beta}$ norms are uniformly bounded.
	Then for $f \equiv 0$ and any $g \in C(\overline{\Omega})$, there exists a unique classical solution $u \in C(\overline{\Omega}) \cap C^2(\Omega)$ to \eqref{def:PDE} with
	\begin{align}\label{def:HJB-classical-solution}
		F(x,r,p,M) \coloneqq \sup_{\alpha \in \mathcal{A}} L^\alpha(x,r,p,M)
	\end{align}
	for $L^\alpha(x,r,p,M) \coloneqq -A^\alpha(x) : M + b^\alpha(x) \cdot p + c^\alpha(x) r + \xi^\alpha(x)$.
	In addition, \eqref{ineq:minimization} holds.
\end{proposition}
\begin{proof}
	The existence of classical solutions to \eqref{def:PDE} can be found in \cite{Safonov1988,GilbargTrudinger2001}.
	The proof of \eqref{ineq:minimization} follows the proof of \Cref{prop:HJB-strong-solution}; further details are therefore omitted.
\end{proof}

We now show that, in the setting of \Cref{prop:HJB-classical-solution}, $f \in L^n(\Omega)$ is allowed.
We claim that (i) a unique generalized solution $u$ to \eqref{def:PDE} exists with $\|u - v\|_{L^\infty(\Omega)} \leq \Psi(v)$ for any $v \in W^{2,n}(\Omega)$ and (ii) \eqref{ineq:minimization} holds.

\paragraph{Ad (i)} This follows from a simple approximation argument.
Let $(f_j)_j \subset C^\infty(\overline{\Omega})$ approximate $f$ with $\lim_{j \to \infty} \|f - f_j\|_{L^n(\Omega)} = 0$. For each $j$, there exists a classical solution $u_j \in C(\overline{\Omega}) \cap W^{2,n}_\mathrm{loc}(\Omega)$ 
to $F[u_j] = f_j$ in $\Omega$ and $u_j = g$ on $\partial \Omega$
(replace $\xi^\alpha$ by $\xi^\alpha-f_j$ in \Cref{prop:HJB-classical-solution}).
\Cref{thm:stability} implies $\|u_j - u_\ell\|_{L^\infty(\Omega)} \leq C\|f_j - f_\ell\|_{L^n(\Omega)}$. Thus, $(u_j)_j$ is a Cauchy sequence in $C(\overline{\Omega})$ with respect to the maximum norm with the limit $u \in C(\overline{\Omega})$. By a similar argument, this limit, called henceforth viscosity solution, is unique in the sense that it does not depend on the choice of $(f_j)_j$.
Since \Cref{thm:stability} holds for classical solutions, it also holds for viscosity solutions as limits of classical solutions. In fact,
given $v \in W^{2,n}(\Omega)$, \Cref{thm:stability} shows
\begin{align*}
	\|u - v\|_{L^\infty(\Omega)} &= \lim_{j \to \infty} \|u_j - v\|_{L^\infty(\Omega)}\\
	&\lesssim \lim_{j \to \infty} (\|g - v\|_{L^\infty(\partial\Omega)} + C\|F[u_j] - F[v]\|_{L^n(\Omega)}) = \Psi(v).
\end{align*}

\paragraph{Ad (ii)}
Let $f \in L^n(\Omega)$ and $\varepsilon>0$ be given.
Then we choose $\widetilde f \in C^\infty(\overline{\Omega})$
such that
$
  \|f - \widetilde{f}\|_{L^n(\Omega)}< \varepsilon/C.
$
Thanks to \Cref{prop:HJB-classical-solution} with 
$f$ replaced by $\widetilde{f}$
(and writing the right-hand side $f$ under the supremum),
the identity 
$
  \inf_{v\in W^{2,n}(\Omega)} \widetilde \Psi(v)=0
$
holds (in analogy to \eqref{ineq:minimization}) for the modified residual
$$
\widetilde\Psi(v)\coloneqq 
	\|g - v\|_{L^\infty(\partial \Omega)} + C\|\widetilde f - F[v]\|_{L^n(\Omega)}.
$$
This and the triangle inequality imply
$$
\inf_{v\in W^{2,n}(\Omega)} \Psi(v)
\leq
\|f-\widetilde f\|_{L^n(\Omega)}
+
\inf_{v\in W^{2,n}(\Omega)} \widetilde \Psi(v)
<
\varepsilon
.
$$
Hence, \eqref{ineq:minimization} holds.

These two points (i)--(ii) are sufficient for convergence of the minimal residual methods from \Cref{sec:FEM}--\ref{sec:NC-FEM} to $u$, defined as the unique limit of classical solutions.

\subsection{Monge--Amp\`ere equation}\label{sub:MA}
In \cite{GallistlTran2023,GallistlTran2024}, a regularized scheme is proposed for the Monge--Amp\`ere equation $\det \D^2 u = \xi$ in $\Omega$ with a right-hand side $\xi \in L^1(\Omega)$.
Given $\varepsilon > 0$, define $\mathbb{S}(\varepsilon) \coloneqq \{A \in \mathbb{S} : \mathrm{tr}(A) = 1 \text{ and } A \geq \varepsilon \mathrm{I}\}$ and
\begin{align}\label{def:HJB-MA}
	F_\varepsilon(x,M) \coloneqq \sup_{A \in \mathbb{S}(\varepsilon)} (-A:M + n\sqrt[n]{\xi(x)\det A}).
\end{align}
If $\xi \in C^{0,\alpha}(\Omega)$, then $F_\varepsilon$ satisfies the assumptions of \Cref{prop:HJB-classical-solution}.
Similar to the previous example, it is possible to define a viscosity solution to \eqref{def:PDE} with $F \coloneqq F_\varepsilon$ from \eqref{def:HJB-MA} by smooth approximation of $\xi \in L^1(\Omega)$ \cite{GallistlTran2024}.

\begin{proposition}[regularized Monge--Amp\`ere equation]
	For any $\xi \in L^n(\Omega)$ and $g \in C(\partial \Omega)$, there exists a unique viscosity solution $u_\varepsilon$ to \eqref{def:PDE} with $F \coloneqq F_\varepsilon$ from \eqref{def:HJB-MA} and $f \equiv 0$. Furthermore, $\|u_\varepsilon - v\|_{L^\infty(\Omega)} \leq \Psi(v)$ for any $v \in W^{2,n}(\Omega)$ and \eqref{ineq:minimization} holds.
\end{proposition}
\begin{proof}
	The assertion follows from the arguments presented in \Cref{sub:HJB-classical}; further details are omitted.
\end{proof}
It is important to note that, as an effect of the degenerate ellipticity
of the Monge--Amp\`ere equation, the constant $C$ in
\eqref{def:Psi} depends on $\varepsilon$ and, hence, is not uniform
with respect to the limit $\varepsilon\to0$.

For convex domains $\Omega$, nonnegative $\xi \geq 0$, and $f \equiv 0$, the sequence $(u_\varepsilon)_\varepsilon$ of viscosity solutions to $F_\varepsilon(x, \D^2 u_\varepsilon) = 0$ in $\Omega$ and $u_\varepsilon = g$ on $\partial \Omega$
converges uniformly to the Alexandrov solution $u$ to the Monge--Amp\`ere equation 
$\det \D^2 u = \xi$
in $\Omega$ and $u = g$ on $\partial \Omega$ if it exists \cite{GallistlTran2023,GallistlTran2024}. We refer to the monograph \cite{Figalli2017} for the precise definition and sufficient conditions for the existence of Alexandrov solutions.

\subsection{Isaac's equation}
This example appeared in \cite{CaffarelliCrandallKocanSwiech1996}.
For countable index sets $\mathcal{A}$ and $\mathcal{B}$, let $A^{\alpha, \beta} \in W^{1,\infty}(\Omega;\mathbb{S})$ with $\lambda \mathrm{I} \leq A^{\alpha, \beta} \leq \Lambda \mathrm{I}$ and $b^{\alpha,\beta} \in W^{1,\infty}(\Omega; \mathbb{R}^n)$ be given such that $\|A^{\alpha,\beta}\|_{W^{1,\infty}(\Omega)} + \|b^{\alpha,\beta}\|_{W^{1,\infty}(\Omega)} \leq C$ for any $\alpha \in \mathcal{A}$, $\beta \in \mathcal{B}$ with a constant $C$ independent of $\alpha$, $\beta$. Define
\begin{align}\label{def:Isaac}
	F(x,p,M) \coloneqq \inf_{\alpha \in \mathcal{A}}\sup_{\beta \in \mathcal{B}} (-A^{\alpha, \beta}(x): M + b^{\alpha,\beta}(x) \cdot p)  .
\end{align}
For any $f \in L^n(\Omega)$ and $g \in C(\partial \Omega)$, there exists a unique $L^n$ viscosity solution $u$ to \eqref{def:PDE} with $F$ from \eqref{def:Isaac} \cite[Example 3.12]{CaffarelliCrandallKocanSwiech1996}.
A simple approximation argument carried out in \cite{CaffarelliCrandallKocanSwiech1996} (similar to that in \Cref{sub:HJB-classical}) leads to the bound $\|u - v\|_{L^\infty(\Omega)} \leq \Psi(v)$ for any $v \in W^{2,n}(\Omega)$.
Due to the nonconvexity of $F$, there is no interior $W^{2,n}$ estimate in the current literature. Therefore, an existence theory for strong solutions is not available.
Thus, the verification of \eqref{ineq:minimization} appears to be not feasible in the current state.
In conclusion,
convergence of minimal residual methods cannot be guaranteed a~priori, but may be checked a~posteriori using the error bounds of \Cref{sec:FEM}.

\section{Numerical examples}\label{sec:numerical-examples}
This section presents numerical experiments in two and three space dimensions.
For $n=2$, the experiments are carried out using the BFS finite element in the
setting of \Cref{sub:cFEM}, while the nonconforming method 
of \Cref{sec:NC-FEM} is utilized for the case $n=3$.

\subsection{Preliminaries}\label{sec:implementation}
The subsequent examples concern Hamilton--Jacobi--Bellman equations of the following form.
Given a compact subset $S$ of $\mathbb{R}^{n \times n}$, we consider
\begin{align}\label{def:HJB-exp}
	F[v] \coloneqq \sup_{A \in S} (-A : \D^2 v + f_A) = 0 \quad\text{in } \Omega,
\end{align}
where $f_A$ depends continuously on $A$.

\subsubsection{Iterative algorithm}
The difficulty in minimizing $\Phi_j$ from \eqref{def:Phi-pract} in $\mathcal{A}_j(g_j)$ and $\Phi_j^\mathrm{nc}$ from \eqref{def:Phi-nc-pract} in $\mathcal{A}_\mathrm{nc}(g_j,\mathcal{T}_j)$ is the non-convexity and non-smoothness of the objective functional in general (with an exception in the linear case).
Therefore, sufficient conditions for minimizers cannot be expected.
Nevertheless, the semismoothness of $F$ \cite{SmearsSueli2014} motivates the following policy iteration. For the sake of brevity, we only provide details for the conforming method of \Cref{sub:cFEM} but mention that a similar algorithm is applied to the nonconforming method of \Cref{sec:NC-FEM} as well.

For any $v \in W^{2,n}(\Omega)$, we denote the supremizing coefficient of \eqref{def:HJB-exp} by
$A_v$, i.e.,
\begin{align*}
	F[v] = - A_v : \D^2 v + f_{A_v}.
\end{align*}
Then the minimization of $\Phi_j$ in $\mathcal{A}_j(g_j)$ from \Cref{sub:cFEM} is carried out as follows.

\begin{algorithm}\label{alg}
	Given an initial input $A_0 \in L^\infty(\Omega;\mathbb{R}^{n\times n})$ with $A_0(x) \in S$ for almost every $x \in \Omega$, compute at each iteration $k = 1,2,\dots$ the discrete minimizer $(s_j^{(k)}, u_j^{(k)}) \in \mathcal{A}_j(g_j)$ of
	\begin{align*}
		\Phi_j^{(k)}(t_j,v_j) \coloneqq \sigma t^n_j + \|f_{A_{k-1}} - A_{k-1}: \D^2 v_j\|_{L^n(\Omega)}^n \quad\text{among } (t_j,v_j) \in \mathcal{A}_j(g_j)
	\end{align*}
	and set $A_k \coloneqq A_{u_j^{(k)}}$.
\end{algorithm}
	In contrast to \cite{SmearsSueli2014}, necessary conditions to the minimization of the nonsmooth functional $\Phi_j$ is not given in form of an operator equation but rather as an inclusion of zero in some generalized subgradient. Since $\Phi_j$ is not convex, the theory of subgradient methods does not apply here. This makes the analysis of our algorithm more challenging.

We note that there is currently no theoretical foundation for this algorithm to guarantee the approximation of the discrete minimizer $u_j$.
In the examples below, the iteration of \Cref{alg} is repeated 8 times with the starting coefficient $A_0 = \mathrm{I}_{n \times n}/2$ and the parameter $\sigma = 10^n$. (For the nonconforming method of \Cref{sec:NC-FEM}, $\tau = 1$.)
The minimization of the objective functional $\Phi_j^{(k)}(t_j,v_j)$ is realized in Matlab using the standard routine \texttt{quadprog} in 2d (and \texttt{fmincon} in 3d for the nonconforming method of \Cref{sec:NC-FEM}).

\subsubsection{Adaptive algorithm}
We utilize a localization of the residual $\|F[u_j]\|_{L^n(\Omega)}^n$ plus some penalization of the boundary residual as refinement indicator. Notice that there is no canonical way to choose the boundary penalization because the $L^\infty$ norm is not localizable. For every $T \in \mathcal{T}_j$, we define the refinement indicator
\begin{align}\label{def:eta_loc}
	\eta_j(T) \coloneqq \sigma\sum_{F \in \mathcal{F}(T) \cap \mathcal{F}(\partial \Omega)} h_F^s \|g - u_j\|_{L^n(F)}^n + \|F[u_j]\|^n_{L^n(T)}.
\end{align}
for some parameter $s$. The quality of the adaptive algorithm appears 
to be sensitive to the choice of $s$; 
we refer to the first computational example
in \Cref{ss:numMA2d} below.
Unless stated otherwise, the default value for $s$ is set to 2.
The standard adaptive loop selects at each refinement step a subset $\mathcal{M}_j \subset \mathcal{T}_j$ of minimal cardinality in 
the D\"orfler marking \cite{Doerfler1996}, i.e.,
\begin{align*}
	\sum_{T \in \mathcal{T}_j} \eta_j(T) \leq \theta \sum_{T \in \mathcal{M}_j} \eta_j(T)
	,
\end{align*}
where we choose $\theta=1/3$.

\subsubsection{Displayed quantities}
All displayed errors in the $L^\infty$, $L^n$, or $W^{1,n}$ norms are 
relative errors.
They are depicted in a log-log plot against the numbers of degrees of freedom ($\mathrm{ndof}$).

\subsection{The regularized Monge--Amp\`ere equation}

In this example, we consider the regularization for the Monge--Amp\`ere equation $\det \D^2 u = \xi$ from \Cref{sub:MA}, where $S = \mathbb{S}(\varepsilon)$ and $f_A \coloneqq n \sqrt[n]{\xi \det A}$ in \eqref{def:HJB-exp}.
We mention that we do not introduce $\varepsilon$-dependent weights in
the definition of $\Psi_j$, in spite of the bad scaling of the 
ellipticity constant with respect to $\varepsilon$.
The reason is that the choice \eqref{def:Phi-pract} turns out to
work well in our examples.

\subsubsection{Monge--Amp\`ere in 2d}\label{ss:numMA2d}
We consider the unit square $\Omega=(0,1)^2$ with the exact solution
$$
u(x,y) = -\left( \left(\sin(\pi x)\right)^{-1} +
\left(\sin(\pi y)\right)^{-1} \right)^{-1}
$$
with homogenous boundary data and right-hand side
$$
\xi(x,y) = \frac{4\pi^2\sin(\pi x)^2\sin(\pi y)^2(2-\sin(\pi x)\sin(\pi y))}{(\sin(\pi x) + \sin(\pi y))^4}
$$
from \cite{GallistlTran2023}. The exact solution $u$ belongs to $H^{2-\delta}(\Omega)$ for any $\delta > 0$, but not to $H^2(\Omega)$.
The convergence with respect to $\varepsilon$ has been analyzed and experimentally investigated in \cite{GallistlTran2023}, so the focus is on the convergence with respect to $\mathrm{ndof}$, in particular, in adaptive computations. Therefore, $\varepsilon = 10^{-4}$ is fixed.

In two dimensions, 
algorithms based on the Cordes condition \cite{SmearsSueli2013,SmearsSueli2014}, e.g.~\cite{GallistlTran2023},
are cheaper for approximating the solution to the regularized problem.
(They are monotone and the boundary data are fixed.)
Reliable and efficient a~posteriori error estimators
are known \cite{GallistlSueli2019,KaweckiSmears2022} and the singularities of $u$ at the corners of the domain suggest an adaptive computation.
However, the error estimators that drive these adaptive algorithms
control the $W^{2,2}$ norm, which scales badly with $\varepsilon$
and is very large for our choice of $\varepsilon$ above.
Therefore, a significant pre-asymptotic range is expected,
similar to \cite[Section 5.4]{GallistlTran2024}.
Even though the solution qualitatively belongs to $W^{2,2}(\Omega)$,
it turns out more effective in this case to
allow variable boundary values of the discrete ansatz functions
for error control in the $L^\infty$ norm.

  \pgfplotstableread{
  ndof	hinv	max_error	L2_error	H1_error	H2_error	energy	eta
   1.6000000000000000e+01   7.0710678118654746e-01   6.1406906250869464e-01   6.5459576999406321e-01   7.2656430931404392e-01   9.3005681955719433e-01   5.6758919518968449e+00   5.6882615636909994e+00
   3.6000000000000000e+01   1.4142135623730949e+00   8.6172868857682508e-02   8.5331818193312689e-02   2.7866714433821754e-01   7.8891117975682423e-01   5.4603694939302483e-01   6.0490435430419121e-01
   1.0000000000000000e+02   2.8284271247461898e+00   5.4043462259472715e-02   5.5130980853996647e-02   2.0107378763556233e-01   7.5695427636172496e-01   3.3435259941074652e-01   2.9950714949717977e-01
   3.2400000000000000e+02   5.6568542494923797e+00   3.3667772972124456e-02   3.3986537773617979e-02   1.3997376936575939e-01   7.3151042571357650e-01   2.2162479389194925e-01   1.7576491513551426e-01
   1.1560000000000000e+03   1.1313708498984759e+01   2.3685108517135289e-02   2.3776040137726141e-02   1.0948070363327951e-01   7.2112258736943613e-01   1.5581276827262341e-01   1.1418911880696131e-01
   4.3560000000000000e+03   2.2627416997969519e+01   1.6703076470093746e-02   1.6640596471905307e-02   7.9615890905114270e-02   7.0568766436259744e-01   1.1079344208987313e-01   7.6851137791235855e-02
   1.6900000000000000e+04   4.5254833995939038e+01   1.1849159658476521e-02   1.1584461643037383e-02   6.5225999812819735e-02   7.0588771186806565e-01   7.7906366201682939e-02   5.2577681479409118e-02
   6.6564000000000000e+04   9.0509667991878075e+01   8.2592773326612100e-03   7.8517817325097616e-03   4.6532434383082487e-02   7.0251178262552083e-01   1.1582701771365805e-01   1.0840165869451751e-01
   }\MAplanarUnifsone
	\pgfplotstableread{
	ndof	hinv	max_error	L2_error	H1_error	H2_error	energy	eta
   1.6000000000000000e+01   7.0710678118654746e-01   6.1406906250869464e-01   6.5459576999406321e-01   7.2656430931404392e-01   9.3005681955719433e-01   5.6758919518968449e+00   5.6882615636909994e+00
   3.6000000000000000e+01   1.4142135623730949e+00   8.6172868857682508e-02   8.5331818193312689e-02   2.7866714433821754e-01   7.8891117975682423e-01   5.4603694939302483e-01   5.1143696798588656e-01
   7.2000000000000000e+01   1.4142135623730949e+00   8.3700133186241921e-02   8.2882088575408394e-02   2.5621046984945989e-01   7.8134656744201780e-01   4.7953589483582187e-01   3.9345962546446345e-01
   8.8000000000000000e+01   1.4142135623730949e+00   7.8112312484127505e-02   7.6146851818656200e-02   2.4170494227461201e-01   7.7707126332955090e-01   4.3448430023061102e-01   3.2868660613812245e-01
   1.0000000000000000e+02   2.8284271247461898e+00   5.4043462259473249e-02   5.5130980853990853e-02   2.0107378763555939e-01   7.5695427636172441e-01   3.3435259941074652e-01   2.2695116700722340e-01
   1.6000000000000000e+02   2.8284271247461898e+00   4.4566649559970965e-02   4.5655436624951570e-02   1.6671073161140298e-01   7.4671022800969322e-01   2.8324537507875625e-01   2.0593616943779330e-01
   1.8000000000000000e+02   2.8284271247461898e+00   3.4112738155572010e-02   3.6042086321013689e-02   1.4522377004744155e-01   7.3337635989394601e-01   2.3831959201920849e-01   1.7632021292604549e-01
   2.7600000000000000e+02   2.8284271247461898e+00   2.4556303206440998e-02   2.8448059472130298e-02   1.1688729386885519e-01   7.2341856219307554e-01   1.9379917450230225e-01   1.5583248469288213e-01
   4.2800000000000000e+02   2.8284271247461898e+00   1.8756406063504927e-02   2.2904781388478637e-02   9.4126442674650923e-02   7.1655850475307092e-01   1.6264282075946829e-01   1.3513712339226788e-01
   5.3600000000000000e+02   5.6568542494923797e+00   1.8401608038408152e-02   2.1135730726886581e-02   8.9691583699254521e-02   7.1309105530391881e-01   1.4147013118055657e-01   1.0939498476544532e-01
   7.0400000000000000e+02   5.6568542494923797e+00   1.4597372420249463e-02   1.7011135721084290e-02   7.3489138283776040e-02   7.1135518148353249e-01   1.1876334962970618e-01   9.4842039218052607e-02
   8.6400000000000000e+02   5.6568542494923797e+00   1.3244226596722823e-02   1.4372787711476862e-02   6.6404905461308922e-02   7.0314935105647691e-01   1.0232179234095003e-01   7.8678621398025630e-02
   1.0720000000000000e+03   5.6568542494923797e+00   1.0345107955813499e-02   1.1582984515239233e-02   5.5108274580600893e-02   7.0469682745881002e-01   8.6111023220496422e-02   6.9315620108004017e-02
   1.2920000000000000e+03   5.6568542494923797e+00   8.9168277494371378e-03   9.7650894012332057e-03   4.9073519124335073e-02   6.9150760767593467e-01   7.4543824839682360e-02   6.0171053630043769e-02
   1.6720000000000000e+03   5.6568542494923797e+00   7.0785669758153674e-03   8.0076215133230832e-03   4.1755675116723887e-02   6.9683011315437726e-01   6.2934869762728954e-02   5.2354462397061922e-02
   2.2200000000000000e+03   5.6568542494923797e+00   5.4259204616380168e-03   6.2095407556029322e-03   3.2773186481825180e-02   6.9785144146734135e-01   5.2091867192547187e-02   4.4605542331895863e-02
   2.8520000000000000e+03   5.6568542494923797e+00   4.8257728448493657e-03   5.2582435530294271e-03   2.8912657458447219e-02   6.8648636160684129e-01   4.3509707080187421e-02   3.6326335617025106e-02
   3.4720000000000000e+03   5.6568542494923797e+00   3.8932006849024558e-03   4.2518534484192287e-03   2.4916741185661385e-02   6.9458607279260798e-01   3.5560771448761803e-02   2.9883330392144400e-02
   4.3560000000000000e+03   5.6568542494923797e+00   3.0484900323028679e-03   3.3764257138908815e-03   2.0952928494131884e-02   6.7502097090615720e-01   2.9424557431052652e-02   2.5267642877242844e-02
   5.5160000000000000e+03   1.1313708498984759e+01   2.4362954836214040e-03   2.6603484307398733e-03   1.6324495628259088e-02   6.8267245650475950e-01   2.4428460482839032e-02   2.1201637998380238e-02
   7.1720000000000000e+03   1.1313708498984759e+01   2.0020874289604636e-03   2.1399097656171697e-03   1.3186220884415387e-02   6.6646670082467607e-01   2.0005953762232010e-02   1.7337750994067966e-02
   9.1160000000000000e+03   1.1313708498984759e+01   1.6105683916966165e-03   1.6724809207306246e-03   1.0909446460702384e-02   6.7664806352393869e-01   1.6373315844098759e-02   1.4270667357870782e-02
   1.1636000000000000e+04   1.1313708498984759e+01   1.3661092587537256e-03   1.4250155942899299e-03   9.7506249439806517e-03   6.6533431659953746e-01   1.3674192893619068e-02   1.1859619047599849e-02
   1.4708000000000000e+04   1.1313708498984759e+01   1.0987139402691610e-03   1.1406362252359020e-03   8.1208685569315545e-03   6.7410904433589769e-01   1.1202355121854038e-02   9.7748032260936805e-03
   1.8184000000000000e+04   1.1313708498984759e+01   8.7290086970958703e-04   8.9748600744115020e-04   6.5006990968638741e-03   6.5916333377666980e-01   9.1672827305575629e-03   8.0690825946961170e-03
   2.2892000000000000e+04   1.1313708498984759e+01   7.2098896975861238e-04   7.3363557207739496e-04   5.4983851401584427e-03   6.7002494122831091e-01   7.5172907789922512e-03   6.6027263808504369e-03
   2.8968000000000000e+04   1.1313708498984759e+01   5.5518582905274587e-04   5.6315387049489833e-04   4.3961942024210278e-03   6.5451480565672837e-01   6.1414625015949941e-03   5.4818517437349931e-03
   3.7056000000000000e+04   1.1313708498984759e+01   4.4661841087123067e-04   4.4979223142458091e-04   3.5526019024697338e-03   6.6366381843367339e-01   5.0096023537775907e-03   4.4870999302462815e-03
   4.7300000000000000e+04   1.1313708498984759e+01   3.6194503051725270e-04   3.6354894012452580e-04   2.9880221947714506e-03   6.5165666088646468e-01   4.1223059212771420e-03   3.7054525283458950e-03
   6.0380000000000000e+04   2.2627416997969519e+01   2.9963195235058475e-04   2.9709577413422828e-04   2.3210814753187804e-03   6.6262082962153623e-01   3.3726297782825094e-03   3.0222649968104222e-03
   7.5692000000000000e+04   2.2627416997969519e+01   2.4850034001635342e-04   2.4475556060109667e-04   1.8323368539971056e-03   6.5263810928499044e-01   2.7691322994998669e-03   2.4751811537286812e-03
   9.4424000000000000e+04   2.2627416997969519e+01   2.0893583542063799e-04   2.0318362372135211e-04   1.5460171982779303e-03   6.4375005992258427e-01   2.2745546278335358e-03   2.0209494818627643e-03
   1.1827600000000000e+05   2.2627416997969519e+01   1.6082512660203015e-04   1.5714695613161726e-04   1.2075782298899568e-03   6.5081218340535441e-01   1.8704142327606466e-03   1.6893641045662407e-03
	}\MAplanarAdapstwo
	
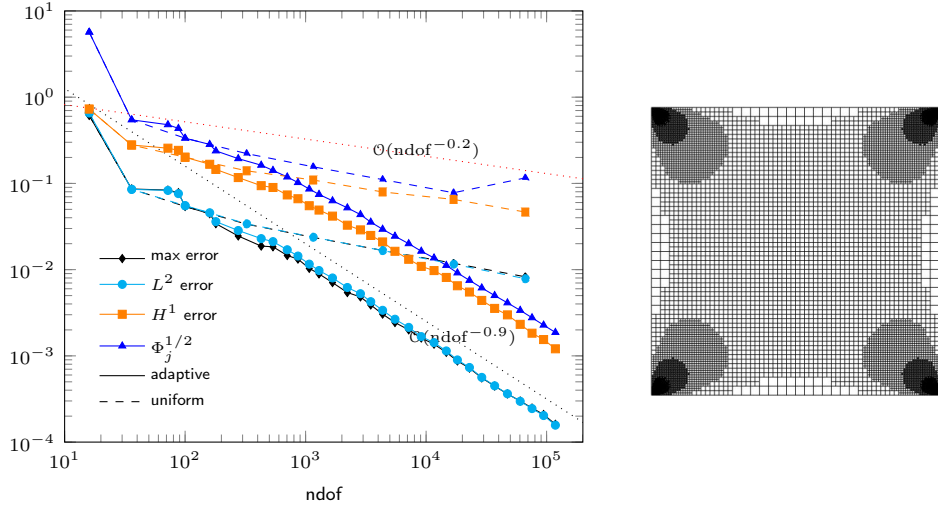
\begin{figure}
	\centering
	\begin{minipage}{.66\textwidth}
\begin{tikzpicture}
		\begin{loglogaxis}[legend pos=south west,legend cell align=left,
			legend style={fill=none,draw=none},
			ymin=1e-4,ymax=1e1,
			xmin=10,xmax=2e5,
			ytick={1e-4,1e-3,1e-2,1e-1,1e0,1e1},
			xtick={1e1,1e2,1e3,1e4,1e5,1e6,1e7}
			]
			\pgfplotsset{
				cycle list={%
					{black, mark=diamond*, mark size=1.5pt},
					{cyan, mark=*, mark size=1.5pt},
					{orange, mark=square*, mark size=1.5pt},
					{blue, mark=triangle*, mark size=1.5pt},
					{black, mark=none, mark size=1.5pt},
					{black, dashed, mark=none, mark size=1.5pt},
					{black, dashed, mark=diamond*, mark size=1.5pt},
					{cyan, dashed, mark=*, mark size=1.5pt},
					{orange, dashed, mark=square*, mark size=1.5pt},
					{blue, dashed, mark=triangle*, mark size=1.5pt},
				},
				legend style={
					at={(0.05,0.05)}, anchor=south west},
				font=\sffamily\scriptsize,
				xlabel=ndof,ylabel=
			}
			\addplot+ table[x=ndof,y=max_error]{\MAplanarAdapstwo};
			\addlegendentry{\tiny max error}
			\addplot+ table[x=ndof,y=L2_error]{\MAplanarAdapstwo};
			\addlegendentry{\tiny $L^2$ error}
			\addplot+ table[x=ndof,y=H1_error]{\MAplanarAdapstwo};
			\addlegendentry{\tiny $H^1$ error}
			\addplot+ table[x=ndof,y=energy]{\MAplanarAdapstwo};
			\addlegendentry{\tiny $\Phi_j^{1/2}$}
			\addplot+ coordinates{(1e0,1e-20) (1e11,1e-20)};
			\addlegendentry{\tiny adaptive}
			\addplot+ coordinates{(1e2,1e-21) (1e11,1e-21)};
			\addlegendentry{\tiny uniform}
			
			\addplot+ table[x=ndof,y=max_error]{\MAplanarUnifsone};
			\addplot+ table[x=ndof,y=L2_error]{\MAplanarUnifsone};
			\addplot+ table[x=ndof,y=H1_error]{\MAplanarUnifsone};
			\addplot+ table[x=ndof,y=energy]{\MAplanarUnifsone};
			
			\addplot+ [dotted,mark=none,black] coordinates{(1e0,1e1) (1e10,1e-8)};
			\node(z) at  (axis cs:2e4,1e-3)
			[above] {\tiny $\mathcal O (\mathrm{ndof}^{-0.9})$};
			\addplot+ [dotted,mark=none,red] coordinates{(1e0,1.3e0) (1e10,1.3e-2)};
			\node(z) at  (axis cs:1e4,1.5e-1)
			[above] {\tiny $\mathcal O (\mathrm{ndof}^{-0.2})$};
		\end{loglogaxis}
	\end{tikzpicture}
	\end{minipage}
	\hfil
	\begin{minipage}{.3\textwidth}
	\includegraphics[width=\textwidth]{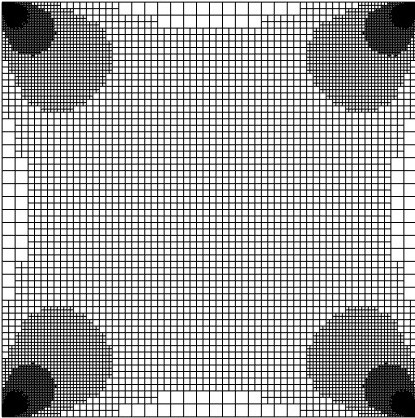}
	\end{minipage}
	\caption{%
	Two-dimensional Monge--Amp\`ere experiment 
	         of \Cref{ss:numMA2d}.
	Convergence history (left) and adaptive mesh (right).}
	\label{fig:ma_exp1}
\end{figure}

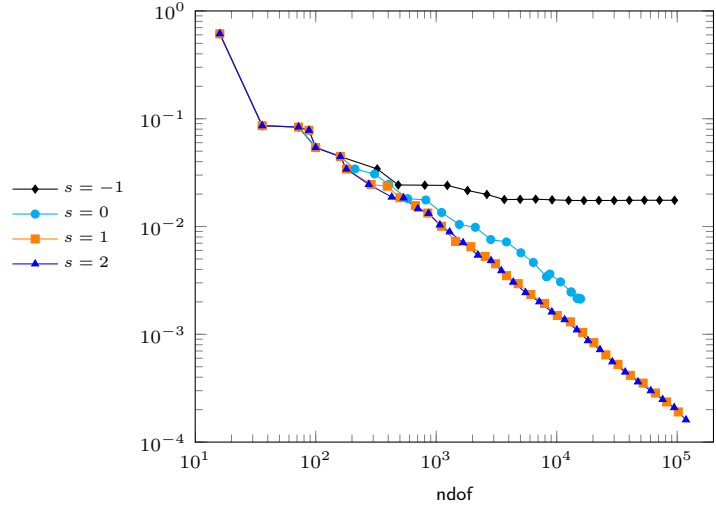
\begin{figure}
\pgfplotstableread{
ndof	hinv	max_error	L2_error	H1_error	H2_error	energy	eta
      1.6000000000000000e+01   7.0710678118654746e-01   6.1406906250869464e-01   6.5459576999406321e-01   7.2656430931404392e-01   9.3005681955719433e-01   5.6758919518968449e+00   5.6882615636909994e+00
   3.6000000000000000e+01   1.4142135623730949e+00   8.6172868857682508e-02   8.5331818193312689e-02   2.7866714433821754e-01   7.8891117975682423e-01   5.4603694939302483e-01   9.9597103956772404e-01
   7.2000000000000000e+01   1.4142135623730949e+00   8.3700133186241921e-02   8.2882088575408394e-02   2.5621046984945989e-01   7.8134656744201780e-01   4.7953589483582187e-01   1.0103356471916687e+00
   1.0000000000000000e+02   2.8284271247461898e+00   5.4043462259473249e-02   5.5130980853990853e-02   2.0107378763555939e-01   7.5695427636172441e-01   3.3435259941074652e-01   9.2393643746723342e-01
   1.6000000000000000e+02   2.8284271247461898e+00   4.4566649559970965e-02   4.5655436624951570e-02   1.6671073161140298e-01   7.4671022800969322e-01   2.8324537507875625e-01   8.2668271674075233e-01
   3.2400000000000000e+02   2.8284271247461898e+00   3.4225024542521881e-02   3.5076235755567058e-02   1.4389421858349172e-01   7.3254770339846909e-01   2.3249754335227990e-01   8.6514219875584109e-01
   4.8400000000000000e+02   2.8284271247461898e+00   2.4296791522797193e-02   2.5738272256086107e-02   1.1017985174663029e-01   7.1901024075989772e-01   1.8008217170484120e-01   7.2395550957334176e-01
   8.0400000000000000e+02   2.8284271247461898e+00   2.4175435189912289e-02   2.4411907638024196e-02   1.0707108245979439e-01   7.1632989292510563e-01   1.7294972776159803e-01   8.7996212058748702e-01
   1.2360000000000000e+03   2.8284271247461898e+00   2.4064325505517667e-02   2.4168750675634995e-02   1.0643076140610296e-01   7.1465110664224607e-01   1.7142687577032953e-01   1.0558632920144888e+00
   1.8200000000000000e+03   2.8284271247461898e+00   2.1591823588957223e-02   2.1728998313529469e-02   9.7234322101235338e-02   7.1469860143149089e-01   1.5914207734674940e-01   1.1361979070781010e+00
   2.6320000000000000e+03   2.8284271247461898e+00   1.9840605470574559e-02   2.0008683316662772e-02   9.0189458522327640e-02   7.1113935169525155e-01   1.5082776937020528e-01   1.2481716250912320e+00
   3.6760000000000000e+03   2.8284271247461898e+00   1.7791242641496768e-02   1.7904637746433319e-02   8.2018223455336728e-02   7.0528269490310047e-01   1.3987181135445750e-01   1.3217338380333314e+00
   4.9720000000000000e+03   2.8284271247461898e+00   1.7856569664751458e-02   1.7925624607476799e-02   8.2460916577306351e-02   7.0578306566980009e-01   1.4153004631530317e-01   1.5478677088861923e+00
   6.6920000000000000e+03   2.8284271247461898e+00   1.7925799120885738e-02   1.7956889608929719e-02   8.2492300017179720e-02   7.0518033548677106e-01   1.4222435486066254e-01   1.7947278820010155e+00
   9.1560000000000000e+03   2.8284271247461898e+00   1.7691657870591990e-02   1.7512654647267350e-02   8.1711864738352777e-02   7.0409974870993419e-01   1.4706936693757913e-01   2.0896917136886319e+00
   1.2564000000000000e+04   2.8284271247461898e+00   1.7479160320228060e-02   1.6981093845825943e-02   7.9472344677735979e-02   7.0049832934791068e-01   1.4050895274731107e-01   2.3957179175561274e+00
   1.6824000000000000e+04   2.8284271247461898e+00   1.7430802934707744e-02   1.6827463267086520e-02   7.9253095345489674e-02   7.0004080354400389e-01   1.4595952723009023e-01   2.7773451539631657e+00
   2.2384000000000000e+04   2.8284271247461898e+00   1.7452035697259127e-02   1.6839864441289917e-02   7.9727211669854539e-02   7.0104232107818254e-01   1.4498548549356904e-01   3.2296562015136243e+00
   2.9944000000000000e+04   2.8284271247461898e+00   1.7448214690116481e-02   1.6833325934650860e-02   7.9734714689444139e-02   7.0081426756848253e-01   1.4187015988972429e-01   3.7412426260588383e+00
   4.0040000000000000e+04   2.8284271247461898e+00   1.7502137038559167e-02   1.6846851922873354e-02   7.9832664215827562e-02   7.0153850656951811e-01   1.4822651889057420e-01   4.3402084594617589e+00
   5.3448000000000000e+04   2.8284271247461898e+00   1.7536530553428634e-02   1.6841963250268970e-02   7.9947670320650796e-02   7.0263667956807463e-01   1.7929104860168477e-01   5.0293004472611766e+00
   7.1304000000000000e+04   2.8284271247461898e+00   1.7519986520281711e-02   1.6774787027029548e-02   7.9821611312795387e-02   7.0206018429274586e-01   1.5105072092012117e-01   5.8126648112739305e+00
   9.5112000000000000e+04   2.8284271247461898e+00   1.7521444545062617e-02   1.6741479110893319e-02   8.0001781318758403e-02   7.0311457200830940e-01   2.2908687422335081e-01   6.7161864154241879e+00
}\MAplanarAdapsminusone
\pgfplotstableread{
ndof	hinv	max_error	L2_error	H1_error	H2_error	energy	eta
   1.6000000000000000e+01   7.0710678118654746e-01   6.1406906250869464e-01   6.5459576999406321e-01   7.2656430931404392e-01   9.3005681955719433e-01   5.6758919518968449e+00   5.6882615636909994e+00
   3.6000000000000000e+01   1.4142135623730949e+00   8.6172868857682508e-02   8.5331818193312689e-02   2.7866714433821754e-01   7.8891117975682423e-01   5.4603694939302483e-01   7.5801866014189778e-01
   7.2000000000000000e+01   1.4142135623730949e+00   8.3700133186241921e-02   8.2882088575408394e-02   2.5621046984945989e-01   7.8134656744201780e-01   4.7953589483582187e-01   6.5255498834528403e-01
   1.0000000000000000e+02   2.8284271247461898e+00   5.4043462259473249e-02   5.5130980853990853e-02   2.0107378763555939e-01   7.5695427636172441e-01   3.3435259941074652e-01   4.9243815264037955e-01
   1.6000000000000000e+02   2.8284271247461898e+00   4.4566649559970965e-02   4.5655436624951570e-02   1.6671073161140298e-01   7.4671022800969322e-01   2.8324537507875625e-01   3.9361483365438449e-01
   2.1200000000000000e+02   2.8284271247461898e+00   3.4191202362668038e-02   3.5644045649076264e-02   1.4389869488714557e-01   7.3273501219425352e-01   2.3641055460983385e-01   3.3401961416779685e-01
   3.0800000000000000e+02   2.8284271247461898e+00   3.0707150303921585e-02   3.2246061887179663e-02   1.3175765405618442e-01   7.3155851449134157e-01   2.1505826509542236e-01   3.0194021173515689e-01
   4.0400000000000000e+02   2.8284271247461898e+00   2.4595102114919505e-02   2.5676688071711917e-02   1.0909104071328500e-01   7.1864031225615888e-01   1.7938594389033516e-01   2.4381862426271123e-01
   5.8000000000000000e+02   2.8284271247461898e+00   1.8141923215987114e-02   1.9342993123113252e-02   8.4915325765819116e-02   7.0930799059235827e-01   1.4625533966089127e-01   1.9095905939352956e-01
   8.2000000000000000e+02   5.6568542494923797e+00   1.7638318611496277e-02   1.8318435643474636e-02   8.2884268160481861e-02   7.0732451671250318e-01   1.2823269778321950e-01   1.7557158572282305e-01
   1.1080000000000000e+03   5.6568542494923797e+00   1.3504440661526767e-02   1.4405527778052298e-02   6.6685902581995610e-02   7.0454055362105583e-01   1.0839189184879562e-01   1.4229518326403232e-01
   1.5560000000000000e+03   5.6568542494923797e+00   1.0441481621148611e-02   1.1110400191656851e-02   5.3469277893013360e-02   7.0418414923220696e-01   8.7592303977082966e-02   1.1343825182045619e-01
   2.1120000000000000e+03   5.6568542494923797e+00   9.8359153423471114e-03   1.0250759393791836e-02   5.0740930514461637e-02   6.9826205368184147e-01   7.8921240548911573e-02   1.0424576650532695e-01
   2.8360000000000000e+03   5.6568542494923797e+00   7.5722443895083402e-03   7.9711026077840599e-03   4.1380187232774548e-02   7.0106552103285291e-01   6.7048323529060636e-02   8.7288345104064247e-02
   3.8240000000000000e+03   5.6568542494923797e+00   7.1980991194364680e-03   7.6061446854041948e-03   4.1342541899542909e-02   6.9915170338543553e-01   6.1251675867849312e-02   8.2976855754294615e-02
   5.0400000000000000e+03   5.6568542494923797e+00   5.7053520652714578e-03   5.8722984133270945e-03   2.9896846863229406e-02   6.9981679457773116e-01   5.1245111228087982e-02   6.1751297805195983e-02
   6.3960000000000000e+03   5.6568542494923797e+00   4.6392515337594573e-03   4.7907788367508196e-03   2.6190012432118309e-02   6.8226649748792745e-01   3.9147478814897455e-02   4.9570517491137479e-02
   8.2280000000000000e+03   1.1313708498984759e+01   3.4222615317049949e-03   3.9984456161228721e-03   2.3897046545572505e-02   6.9985816288689218e-01   7.0389213748952587e-02   7.2714354833404146e-02
   8.3120000000000000e+03   1.1313708498984759e+01   3.4421606109023175e-03   3.9838727769323130e-03   2.2023854050888962e-02   6.9813305164738415e-01   4.8577658318250755e-02   5.1348458847387250e-02
   8.7840000000000000e+03   1.1313708498984759e+01   3.6269755783761193e-03   3.8395161892015427e-03   2.3327706436960135e-02   6.9775232830208422e-01   3.6146901125086665e-02   4.4289459917123700e-02
   1.0752000000000000e+04   1.1313708498984759e+01   3.0661949419961499e-03   3.2194214813424625e-03   2.0126001015606354e-02   6.7922142504859095e-01   2.9593024875643928e-02   3.8033692109503740e-02
   1.3168000000000000e+04   1.1313708498984759e+01   2.4704452011812691e-03   2.6572387329660331e-03   1.6250925617626307e-02   6.8904748422563278e-01   6.1589585481450868e-02   6.3488623847936992e-02
   1.3200000000000000e+04   1.1313708498984759e+01   2.4692867177136596e-03   2.6471033196259880e-03   1.5124213442967284e-02   6.8617955086532390e-01   2.8767644398039287e-02   3.2172843150404903e-02
   1.4764000000000000e+04   1.1313708498984759e+01   2.1568400598499450e-03   2.3888469057669597e-03   1.7255743568110726e-02   6.8747850879421080e-01   1.6277351169966500e-01   1.6352779201153256e-01
   1.4864000000000000e+04   1.1313708498984759e+01   2.1640702377215991e-03   2.3983980529249851e-03   1.7259314188464492e-02   6.8733222174788078e-01   1.3578987513141536e-01   1.3670920614967561e-01
   1.5016000000000000e+04   1.1313708498984759e+01   2.1605297838582552e-03   2.3927058288425156e-03   1.7326132544446652e-02   6.8816976484554704e-01   1.3404527067587074e-01   1.3497403030705102e-01
   1.5176000000000000e+04   1.1313708498984759e+01   2.1537190152969691e-03   2.3999096177651645e-03   1.7553602332249733e-02   6.8858066521251693e-01   1.6031273470469884e-01   1.6106908142685525e-01
   1.5416000000000000e+04   1.1313708498984759e+01   2.1423108660820508e-03   2.3296221676321082e-03   1.5904405002045111e-02   6.8418557480485143e-01   7.2354860411989533e-02   7.4067456349666175e-02
   1.5500000000000000e+04   1.1313708498984759e+01   2.1414278876521592e-03   2.3295060692581308e-03   1.5379933906638953e-02   6.8303661554309669e-01   3.5307086150966602e-02   3.8308361747529426e-02
   1.5696000000000000e+04   1.1313708498984759e+01   2.1281607723994699e-03   2.3692243015317306e-03   1.5934238733265653e-02   6.8453322024200203e-01   1.0997227445089379e-01   1.1092062860101044e-01
   1.5776000000000000e+04   1.1313708498984759e+01   2.1371607096553271e-03   2.3840730967826904e-03   1.6480614630061910e-02   6.8508739635929683e-01   9.6860943259412052e-02   9.8119301537607659e-02
}\MAplanarAdapszero
\pgfplotstableread{
ndof	hinv	max_error	L2_error	H1_error	H2_error	energy	eta
   1.6000000000000000e+01   7.0710678118654746e-01   6.1406906250869464e-01   6.5459576999406321e-01   7.2656430931404392e-01   9.3005681955719433e-01   5.6758919518968449e+00   5.6882615636909994e+00
   3.6000000000000000e+01   1.4142135623730949e+00   8.6172868857682508e-02   8.5331818193312689e-02   2.7866714433821754e-01   7.8891117975682423e-01   5.4603694939302483e-01   6.0490435430419121e-01
   7.2000000000000000e+01   1.4142135623730949e+00   8.3700133186241921e-02   8.2882088575408394e-02   2.5621046984945989e-01   7.8134656744201780e-01   4.7953589483582187e-01   4.7797951126733718e-01
   8.8000000000000000e+01   1.4142135623730949e+00   7.8112312484127505e-02   7.6146851818656200e-02   2.4170494227461201e-01   7.7707126332955090e-01   4.3448430023061102e-01   4.1196302967673110e-01
   1.0000000000000000e+02   2.8284271247461898e+00   5.4043462259473249e-02   5.5130980853990853e-02   2.0107378763555939e-01   7.5695427636172441e-01   3.3435259941074652e-01   2.9950714949718049e-01
   1.6000000000000000e+02   2.8284271247461898e+00   4.4566649559970965e-02   4.5655436624951570e-02   1.6671073161140298e-01   7.4671022800969322e-01   2.8324537507875625e-01   2.4715796194092757e-01
   1.8000000000000000e+02   2.8284271247461898e+00   3.4112738155572010e-02   3.6042086321013689e-02   1.4522377004744155e-01   7.3337635989394601e-01   2.3831959201920849e-01   2.0966368352999529e-01
   2.8800000000000000e+02   2.8284271247461898e+00   2.4578426529852235e-02   2.8100786472385557e-02   1.1589308708398287e-01   7.2281276529861560e-01   1.9203296278356718e-01   1.7143828233139147e-01
   3.9200000000000000e+02   2.8284271247461898e+00   2.3682032982867578e-02   2.5468675896441171e-02   1.0705355293103508e-01   7.2050634696870752e-01   1.7679145239229532e-01   1.4792785699381195e-01
   5.0000000000000000e+02   5.6568542494923797e+00   1.8500350076986326e-02   2.1319990034010701e-02   8.9963780058825438e-02   7.1337122721363266e-01   1.4305059500938591e-01   1.2073748358678955e-01
   6.7200000000000000e+02   5.6568542494923797e+00   1.5609222438606674e-02   1.7551258850994622e-02   7.5869094403750767e-02   7.1212936165880369e-01   1.2360259552323496e-01   1.0329795517336393e-01
   8.5200000000000000e+02   5.6568542494923797e+00   1.3336083601176376e-02   1.4466357837709803e-02   6.6732264154865814e-02   7.0368935328970372e-01   1.0330772012285304e-01   8.5185540762712772e-02
   1.1040000000000000e+03   5.6568542494923797e+00   1.0071766878703915e-02   1.1152905728973986e-02   5.3915974227430684e-02   7.0243017045831224e-01   8.3677845239866058e-02   7.1195138421463494e-02
   1.4480000000000000e+03   5.6568542494923797e+00   7.2790219878828338e-03   8.3802093958528501e-03   4.2836262273353326e-02   6.9885539091748161e-01   6.7623458828335634e-02   5.9974461005668250e-02
   1.9400000000000000e+03   5.6568542494923797e+00   6.4897811860050161e-03   7.2335593188346467e-03   3.7341762956130974e-02   6.8813180753676539e-01   5.7799201815045699e-02   4.9668938077995085e-02
   2.5480000000000000e+03   5.6568542494923797e+00   5.2777664388282720e-03   5.7062926903567309e-03   3.0718017240284425e-02   6.9436794587238404e-01   4.8117901119676787e-02   4.1667929566722044e-02
   3.1120000000000000e+03   5.6568542494923797e+00   4.5112112415626208e-03   4.9208150820474620e-03   2.7060089939345142e-02   6.7992150146230079e-01   4.0150384295323777e-02   3.4348421990648997e-02
   3.8320000000000000e+03   1.1313708498984759e+01   3.5034090931520229e-03   3.8283208815372526e-03   2.1745909406671121e-02   6.8521059531311956e-01   3.3481889825553432e-02   2.9241628918444142e-02
   4.7840000000000000e+03   1.1313708498984759e+01   2.9543260251567380e-03   3.1559533463969510e-03   1.8773248384904569e-02   6.7220785114712400e-01   2.7667254103686853e-02   2.3989354515865292e-02
   6.1040000000000000e+03   1.1313708498984759e+01   2.3461137781105685e-03   2.5157692493416035e-03   1.5107748617546629e-02   6.7982811834307133e-01   2.3038190790569457e-02   2.0188045697589804e-02
   7.9240000000000000e+03   1.1313708498984759e+01   1.9349460878070569e-03   2.0096090202103957e-03   1.2619394971013904e-02   6.6478800565634211e-01   1.8847310304520675e-02   1.6495036835592662e-02
   1.0100000000000000e+04   1.1313708498984759e+01   1.4969344704983659e-03   1.5626802262142931e-03   1.0127607032027865e-02   6.7125770698669163e-01   1.5359696153230950e-02   1.3634529760328152e-02
   1.2992000000000000e+04   1.1313708498984759e+01   1.3049917467050148e-03   1.3446642695336083e-03   9.0119958829274049e-03   6.6229123973269621e-01   1.2647584630595793e-02   1.1030523096538227e-02
   1.6412000000000000e+04   1.1313708498984759e+01   1.0427081433877316e-03   1.0482904841476876e-03   7.1060078958369403e-03   6.7109837003507966e-01   1.0256679565095676e-02   8.9679565799316539e-03
   2.0348000000000000e+04   2.2627416997969519e+01   8.3620502797000231e-04   8.4405155984316636e-04   5.9977266545831239e-03   6.5667793542485242e-01   8.4691984343527888e-03   7.4515527805671632e-03
   2.5576000000000000e+04   2.2627416997969519e+01   6.4458208952971794e-04   6.4693115833489891e-04   4.6112865517903981e-03   6.6349233907486982e-01   6.8919630179747849e-03   6.1440501329571070e-03
   3.2260000000000000e+04   2.2627416997969519e+01   5.2574972168131643e-04   5.2521221219887799e-04   3.8650483981321395e-03   6.5123699156208559e-01   5.6758851322503178e-03   5.0732770170779078e-03
   4.0864000000000000e+04   2.2627416997969519e+01   4.1669678648866942e-04   4.1339603310963148e-04   3.0498051147017165e-03   6.5997233312893488e-01   4.6369167446963415e-03   4.1705169082687467e-03
   5.2160000000000000e+04   2.2627416997969519e+01   3.5269390170381942e-04   3.4708639984691548e-04   2.6588161853713086e-03   6.5042186453050488e-01   3.8225389117511432e-03   3.4212142014025766e-03
   6.5972000000000000e+04   2.2627416997969519e+01   2.8543896312102786e-04   2.8094881771457584e-04   2.1865324407558755e-03   6.6016614084342451e-01   3.1324222552446901e-03   2.8090885977566883e-03
   8.1772000000000000e+04   2.2627416997969519e+01   2.3574677334504038e-04   2.2993801660420632e-04   1.5797429229765713e-03   6.4998434179279063e-01   2.5766216777014157e-03   2.2959202980218029e-03
   1.0204000000000000e+05   2.2627416997969519e+01   1.9046815078508869e-04   1.8470926546281969e-04   1.3585828677887587e-03   6.5958517550746820e-01   2.1119942738662151e-03   1.8921526631597088e-03
}\MAplanarAdapsone
	\centering
	\begin{tikzpicture}
		\begin{loglogaxis}[legend pos=south west,legend cell align=left,
			legend style={fill=none,draw=none},
			ymin=1e-4,ymax=1e0,
			xmin=10,xmax=2e5,
			ytick={1e-4,1e-3,1e-2,1e-1,1e0,1e1},
			xtick={1e1,1e2,1e3,1e4,1e5,1e6,1e7}
			]
			\pgfplotsset{
				cycle list={%
					{black, mark=diamond*, mark size=1.5pt},
					{cyan, mark=*, mark size=1.5pt},
					{orange, mark=square*, mark size=1.5pt},
					{blue, mark=triangle*, mark size=1.5pt},
					{black, mark=none, mark size=1.5pt},
					{black, dashed, mark=none, mark size=1.5pt},
					{black, dashed, mark=diamond*, mark size=1.5pt},
					{cyan, dashed, mark=*, mark size=1.5pt},
					{orange, dashed, mark=square*, mark size=1.5pt},
					{blue, dashed, mark=triangle*, mark size=1.5pt},
				},
				legend style={
					at={(-0.11,.5)}, anchor=east},
				font=\sffamily\scriptsize,
				xlabel=ndof,ylabel=
			}
			\addplot+ table[x=ndof,y=max_error]{\MAplanarAdapsminusone};
			\addlegendentry{\tiny $s = -1$}
			\addplot+ table[x=ndof,y=max_error]{\MAplanarAdapszero};
			\addlegendentry{\tiny $s = 0$}
			\addplot+ table[x=ndof,y=max_error]{\MAplanarAdapsone};
			\addlegendentry{\tiny $s = 1$}
			\addplot+ table[x=ndof,y=max_error]{\MAplanarAdapstwo};
			\addlegendentry{\tiny $s = 2$}
		\end{loglogaxis}
	\end{tikzpicture}
	\caption{Convergence history plot of the $L^\infty$ error for 
	         two-dimensional Monge--Amp\`ere experiment 
	         of \Cref{ss:numMA2d}
	         with different parameters $s$}
	\label{fig:ma_exp1_comp}
\end{figure}

The convergence history plot of 
\Cref{fig:ma_exp1} reveals the experimental convergence rate 1/5
for all displayed errors on uniform meshes. 
Adaptive computations refine towards the corner singularities as shown in 
the mesh of \Cref{fig:ma_exp1}
and improve the convergence rate of the errors to almost 1.
We also mention that in this example, $\Phi_j^{1/2}$ appears to be efficient as well. The performance of the adaptive algorithm appears to depend on the choice of $s$ in \eqref{def:eta_loc}, which balances refinements towards the boundary or the corner singularity.
While the best choice of $s$ is not known, heuristically, higher-order discretizations require stronger refinements towards the corner singularities.
This is reflected in the convergence history plot of \Cref{fig:ma_exp1_comp} with different $s$. We expect that, for higher-order methods, $s$ can be set to $\infty$, where the convention to $h_F^\infty = 0$ is utilized.
(Recall the default value $s = 2$ for \Cref{fig:ma_exp1}.)

\begin{figure}
	\pgfplotstableread{
ndof	max_error	L3_error	H1_error	H2_error	energy
   5.0000000000000000e+01   1.0232422492029083e-01   7.8723223465985703e-02   1.9776354968474436e-01   5.1756724410759070e-01   8.3819377426164854e-01
   1.0000000000000000e+02   3.7692435099028310e-02   2.1074948917704769e-02   1.1013193030180693e-01   3.2653229355865343e-01   3.5976242583907614e-01
   2.0000000000000000e+02   3.6441201297379225e-02   8.1602428514509287e-03   5.9340289581087620e-02   2.5998288362894895e-01   2.4493121839855803e-01
   4.0000000000000000e+02   1.6297447124948371e-02   9.6766593728687001e-03   5.9254972224250801e-02   2.7627560606445606e-01   2.2263155172878807e-01
   8.0000000000000000e+02   1.2091746947345220e-02   7.1936557730015701e-03   4.3301703799453070e-02   2.6579408227648149e-01   1.9395208573306358e-01
   1.6000000000000000e+03   1.2763787299822815e-02   7.3552519110494053e-03   3.8448458317693668e-02   2.4597395307067124e-01   1.6792639858008468e-01
   3.2000000000000000e+03   1.2520553420122253e-02   6.4402708370300561e-03   3.1582283669721385e-02   2.3516100916514809e-01   1.5233122364223908e-01
   6.4000000000000000e+03   6.9263076819962134e-03   4.6191956057623339e-03   2.1560946128040683e-02   2.2371779625251598e-01   1.3152045075427934e-01
   1.2800000000000000e+04   7.1429301861125963e-03   4.4940412344216850e-03   1.7221337586178359e-02   1.9118553353913231e-01   1.1209640376721096e-01
   2.5600000000000000e+04   4.6332178845737660e-03   3.8833319753105383e-03   1.4367642898320798e-02   2.0567712378641206e-01   1.0123536217689824e-01
   5.1200000000000000e+04   2.5625822175728351e-03   3.0027632614450612e-03   9.9221527187731504e-03   1.7747113203589776e-01   8.6896352500984098e-02
	}\unif

	\pgfplotstableread{
ndof	max_error	L3_error	H1_error	H2_error	energy
   5.0000000000000000e+01   1.0232422492029083e-01   7.8723223465985703e-02   1.9776354968474436e-01   5.1756724410759070e-01   8.3819377426164854e-01
   8.0000000000000000e+01   6.8661177178580252e-02   3.1145380783561431e-02   1.5659451378869110e-01   4.7526284165325905e-01   6.4750885558693905e-01
   9.0000000000000000e+01   9.4552523234839939e-02   3.3842015843472871e-02   1.4265966831538990e-01   4.9841561953139296e-01   5.5565067771141052e-01
   1.0000000000000000e+02   7.5263001800757748e-02   2.8522068066254857e-02   1.0152010353388113e-01   3.7518805597462446e-01   3.6063224438173913e-01
   1.4000000000000000e+02   3.3890156028213603e-02   9.8930046210548104e-03   5.8298243457193859e-02   2.4767969826203390e-01   2.8499630689380839e-01
   1.6000000000000000e+02   3.3400421667515917e-02   9.0448121007200965e-03   5.7728126951932295e-02   2.5858587659564214e-01   2.5251106553424679e-01
   1.9000000000000000e+02   3.9304667230993587e-02   9.2951866975643366e-03   6.6424932620593680e-02   2.7708420702129777e-01   2.4476810244650432e-01
   2.8000000000000000e+02   1.9851946190264560e-02   1.0991072664166213e-02   5.6707928158910327e-02   2.5607477941499412e-01   2.2973889676598100e-01
   3.6000000000000000e+02   1.3963306816121391e-02   9.7015421052103028e-03   5.3651020730217297e-02   2.3200005058431653e-01   2.1994968069365367e-01
   4.7000000000000000e+02   1.5811218736727802e-02   1.0825552247055458e-02   5.5319555822645444e-02   2.6900856225101555e-01   2.0471226784523208e-01
   6.6000000000000000e+02   1.5608931132458502e-02   9.5950127564796359e-03   4.8268741812619441e-02   2.4732512652263389e-01   1.9089965971192691e-01
   9.0000000000000000e+02   1.0199523992390340e-02   5.8698986991983261e-03   3.6458708298393351e-02   2.2689938527131776e-01   1.7695502511635153e-01
   1.2900000000000000e+03   9.1596034840329132e-03   4.9662910577181769e-03   3.0495047410894326e-02   2.3534378469857276e-01   1.5987270121355024e-01
   2.0000000000000000e+03   7.8786827414072148e-03   4.5119935592282320e-03   2.4921906720284322e-02   2.0965358702519390e-01   1.4209515859814703e-01
   2.6800000000000000e+03   5.7987682882615372e-03   3.2116763840806997e-03   2.0933599503687405e-02   2.2542242946566890e-01   1.2901952172432130e-01
   3.9100000000000000e+03   3.6109840106797958e-03   2.1048213787440658e-03   1.6061531625683152e-02   1.9203239520802387e-01   1.1720855909269400e-01
   5.2900000000000000e+03   2.9286725162677081e-03   1.6056261104156339e-03   1.3767800192603505e-02   1.8050453607931399e-01   1.0504913234734288e-01
   7.2600000000000000e+03   2.4045461402077248e-03   1.6075134579722047e-03   1.1739538458621973e-02   1.6378119650745535e-01   9.4801153675394600e-02
   9.5000000000000000e+03   1.6851801349012654e-03   1.3427695572823170e-03   1.0206398143641318e-02   1.5543175694934355e-01   8.4836434341784145e-02
   1.3110000000000000e+04   1.4298338381441703e-03   9.9096309638112496e-04   8.7046645719002797e-03   1.4276707630312427e-01   7.6553481011378058e-02
   1.7360000000000000e+04   1.3361873071144699e-03   7.9166476228666085e-04   7.3910986214524186e-03   1.2711031237068374e-01   6.9608895894506098e-02
   2.3350000000000000e+04   9.9084973199125995e-04   6.6911373298658076e-04   6.2498206877225643e-03   1.2672203299362239e-01   6.2624072282273097e-02
   3.1050000000000000e+04   8.2758731144388659e-04   5.8138954370230104e-04   5.4708687240082763e-03   1.2887725285936452e-01   5.6579614826599185e-02
	}\adap
	
	\begin{tikzpicture}
		\begin{loglogaxis}[legend pos=south west,legend cell align=left,
			legend style={fill=none,draw=none},
			ymin=1e-4,ymax=1e0,
			xmin=1e1,xmax=1e5,
			ytick={1e-4,1e-3,1e-2,1e-1,1e0,1e1},
			xtick={1e1,1e2,1e3,1e4,1e5,1e6,1e7}
			]
			\pgfplotsset{
				cycle list={%
					{black, mark=diamond*, mark size=1.5pt},
					{black, mark=diamond, mark size=1.5pt},
					{cyan, mark=*, mark size=1.5pt},
					{cyan, mark=o, mark size=1.5pt},
					{orange, mark=square*, mark size=1.5pt},
					{orange, mark=square, mark size=1.5pt},
					{blue, mark=triangle*, mark size=1.5pt},
					{blue, mark=triangle, mark size=1.5pt},					
					{black, mark=none, mark size=1.5pt},
				},
				legend style={
					at={(-0.11,.5)}, anchor=east},
				font=\sffamily\scriptsize,
				xlabel=$\mathsf{ndof}$,ylabel=
			}

			\addplot+ table[x=ndof,y=max_error]{\unif};
			\addlegendentry{max error, uniform}
			\addplot+ table[x=ndof,y=max_error]{\adap};
			\addlegendentry{max error, adaptive}
			\addplot+ table[x=ndof,y=L3_error]{\unif};
			\addlegendentry{$L^3$ error, uniform}
			\addplot+ table[x=ndof,y=L3_error]{\adap};
			\addlegendentry{$L^3$ error, adaptive}
			\addplot+ table[x=ndof,y=H1_error]{\unif};
			\addlegendentry{$W^{1,3}$ error, uniform}
			\addplot+ table[x=ndof,y=H1_error]{\adap};
			\addlegendentry{$W^{1,3}$ error, adaptive}
			\addplot+ table[x=ndof,y=energy]{\unif};
			\addlegendentry{$\Phi_j^{1/3}$, uniform}
			\addplot+ table[x=ndof,y=energy]{\adap};
			\addlegendentry{$\Phi_j^{1/3}$, adaptive}	
			
			\addplot+ [color=black,dotted,mark=none] coordinates{(1e1,5e-1) (1e5,5e-3)};
			\addlegendentry{rate $\mathsf{ndof}^{-1/2}$}
			\addplot+ [color=black,dashed,mark=none] coordinates{(1,3e-1) (1e5,3e-4)};
			\addlegendentry{rate $\mathsf{ndof}^{-2/3}$}			
		\end{loglogaxis}
	\end{tikzpicture}
	\caption{%
		Convergence history plot of piecewise quadratic DG discretization for
		the three-dimensional Monge--Amp\`ere example from
		\Cref{ss:numMA3d}.
	}
	\label{f:convMA3d}
\end{figure}
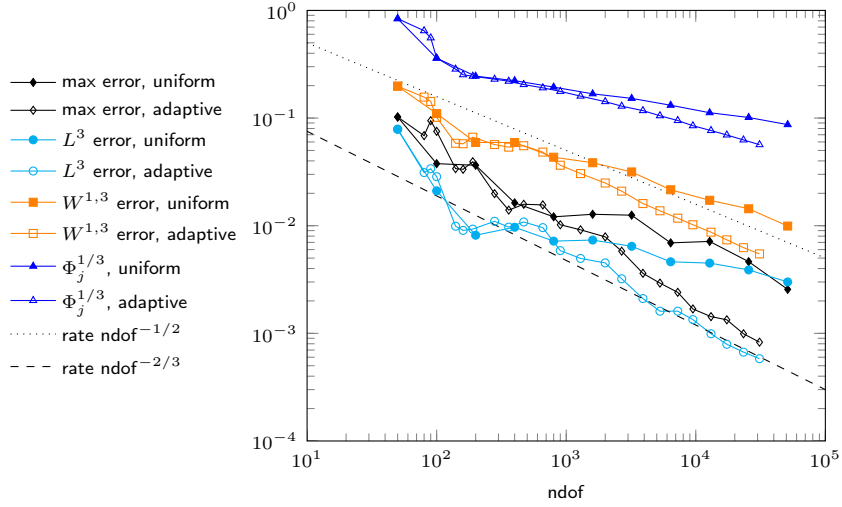

\begin{figure}
	\begin{center}
		\includegraphics[width=.3\textwidth]{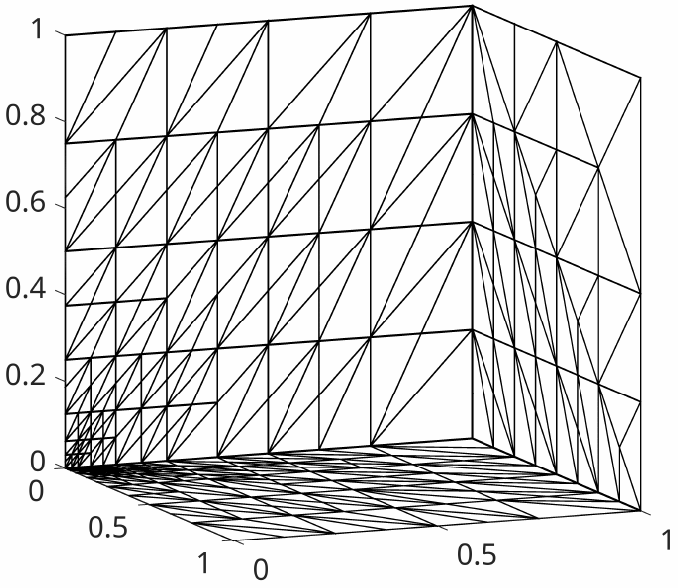}
	\end{center}
	\caption{Adaptive mesh with 3\,105 simplices and 31\,050 degrees of freedom
		in the 23rd iteration of the adaptive algorithm for the
		Monge--Amp\`ere example from \Cref{ss:numMA3d}.
	}
	\label{f:meshMA3d}
\end{figure}

\subsubsection{Monge--Amp\`ere in 3d}\label{ss:numMA3d}
We consider the unit cube $\Omega=(0,1)^3$ 
with $f(x) = 6\cdot(4/5)^4|x|^{-6/5}$
and Dirichlet data such that
the exact solution is given by
$$
  u(x) = |x|^{8/5} .
$$
The solution belongs to the space 
$W^{2+2/5-\nu,3}(\Omega)$ for every $\nu>0$.
The experiment is carried out with a discontinuous Galerkin method
with $k=2$ (piecewise quadratic) 
and the regularization parameter $\varepsilon = 10^{-3}$.
The convergence history is displayed in
\Cref{f:convMA3d}. 
Due to the point singularity, uniform mesh refinement performs
suboptimally and is outperformed by the adaptive algorithm,
which shows a convergence rate close to $2/3$ in terms of
degrees of freedom.
The adaptive process strongly refines towards the singularity
as can be seen from the adaptive mesh displayed in
\Cref{f:meshMA3d}

\subsection{The Pucci equation}\label{ss:numPucci2d}
This example approximates the unknown solution $u \in C(\overline{\Omega}) \cap W^{2,n}_\mathrm{loc}(\Omega)$ to the Pucci equation
\begin{align*}
	\mathcal{P}^+_{\lambda,\Lambda} (\D^2 u) = 1 \text{ in } \Omega \quad\text{and}\quad u = 0 \text{ on } \partial \Omega
\end{align*}
in the L-shaped domain $\Omega = (-1,1)^2 \setminus [0,1]\times[0,-1]$
with the parameters $\lambda = 0.1$ and $\Lambda = 0.9$. Recall the definition of $\mathcal{P}^+_{\lambda,\Lambda}$ from \eqref{def:pucci} to infer $S \coloneqq \{A \in \mathbb{S} : \lambda \mathrm{I}_{2 \times 2} \leq A \leq \Lambda \mathrm{I}_{2 \times 2}\}$ and $f_A \equiv 1$ in \eqref{def:HJB-exp}.
Since the solution is unknown, the objective functional $\Phi_j$ is the only indicator for the convergence of the error.

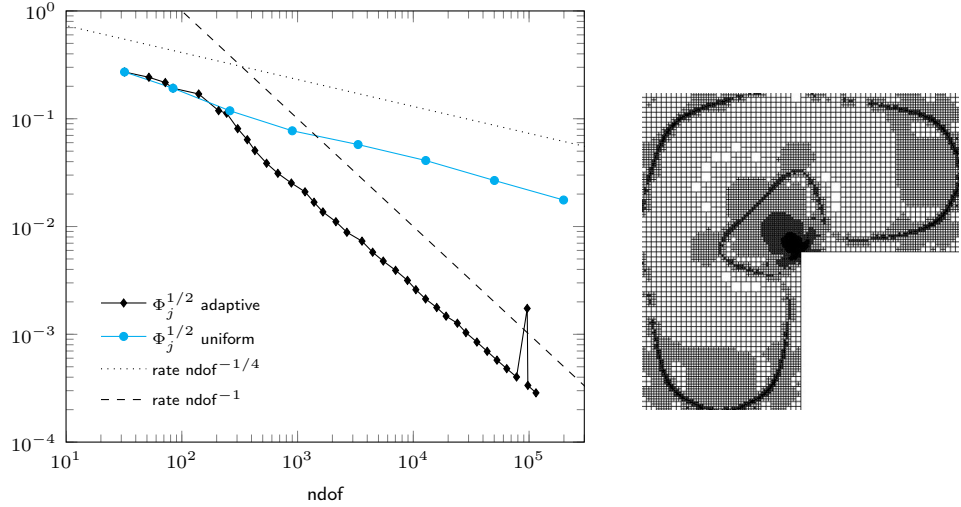
\begin{figure}
	\pgfplotstableread{
ndof	hinv	max_error	L2_error	H1_error	H2_error	energy	eta
   3.2000000000000000e+01   7.0710678118654746e-01                      Inf                      Inf                      Inf                      Inf   2.7143549826977331e-01   5.0746365699524221e-01
   8.4000000000000000e+01   1.4142135623730949e+00                      Inf                      Inf                      Inf                      Inf   1.9193816486150742e-01   1.9298386737042952e-01
   2.6000000000000000e+02   2.8284271247461898e+00                      Inf                      Inf                      Inf                      Inf   1.1855343140375235e-01   8.9982835371778730e-02
   9.0000000000000000e+02   5.6568542494923797e+00                      Inf                      Inf                      Inf                      Inf   7.7291309939243338e-02   5.5903115731989785e-02
   3.3320000000000000e+03   1.1313708498984759e+01                      Inf                      Inf                      Inf                      Inf   5.7606721083197657e-02   4.0718769426287438e-02
   1.2804000000000000e+04   2.2627416997969519e+01                      Inf                      Inf                      Inf                      Inf   4.0961380097290542e-02   2.5605560254101524e-02
   5.0180000000000000e+04   4.5254833995939038e+01                      Inf                      Inf                      Inf                      Inf   2.6735873116077651e-02   1.4730861395213922e-02
   1.9866000000000000e+05   9.0509667991878075e+01                      Inf                      Inf                      Inf                      Inf   1.7617086722859908e-02   1.0343929563543870e-02
    }\pucciunif
   %
	\pgfplotstableread{
ndof	hinv	max_error	L2_error	H1_error	H2_error	energy	eta
   3.2000000000000000e+01   7.0710678118654746e-01                      Inf                      Inf                      Inf                      Inf   2.7143549826977331e-01   5.0746365699524221e-01
   5.2000000000000000e+01   7.0710678118654746e-01                      Inf                      Inf                      Inf                      Inf   2.4241710705764211e-01   3.7469365695901946e-01
   7.2000000000000000e+01   7.0710678118654746e-01                      Inf                      Inf                      Inf                      Inf   2.1662799286475176e-01   2.4638180359480863e-01
   8.4000000000000000e+01   1.4142135623730949e+00                      Inf                      Inf                      Inf                      Inf   1.9193816486150730e-01   1.9298386737042958e-01
   1.4000000000000000e+02   1.4142135623730949e+00                      Inf                      Inf                      Inf                      Inf   1.6967214482991830e-01   1.6776357305380218e-01
   2.0800000000000000e+02   1.4142135623730949e+00                      Inf                      Inf                      Inf                      Inf   1.1885228468028437e-01   1.1282054242786367e-01
   2.4400000000000000e+02   1.4142135623730949e+00                      Inf                      Inf                      Inf                      Inf   1.1244942946195216e-01   1.0447111573917073e-01
   3.0400000000000000e+02   1.4142135623730949e+00                      Inf                      Inf                      Inf                      Inf   8.0895894245643127e-02   7.1693331751431041e-02
   3.6800000000000000e+02   1.4142135623730949e+00                      Inf                      Inf                      Inf                      Inf   6.3889558488446768e-02   5.3843473288402117e-02
   4.2800000000000000e+02   1.4142135623730949e+00                      Inf                      Inf                      Inf                      Inf   5.0673035759010721e-02   4.2111421118123479e-02
   5.4000000000000000e+02   1.4142135623730949e+00                      Inf                      Inf                      Inf                      Inf   3.8608306674137720e-02   3.2543687866575889e-02
   6.7600000000000000e+02   2.8284271247461898e+00                      Inf                      Inf                      Inf                      Inf   3.1147926134763297e-02   2.7681702623878025e-02
   8.8400000000000000e+02   2.8284271247461898e+00                      Inf                      Inf                      Inf                      Inf   2.5325632227444106e-02   2.2526638817455667e-02
   1.1600000000000000e+03   2.8284271247461898e+00                      Inf                      Inf                      Inf                      Inf   2.0983356215076708e-02   1.7878854506013363e-02
   1.3880000000000000e+03   2.8284271247461898e+00                      Inf                      Inf                      Inf                      Inf   1.6791361640661362e-02   1.4783815647866145e-02
   1.6560000000000000e+03   2.8284271247461898e+00                      Inf                      Inf                      Inf                      Inf   1.3654465763883938e-02   1.2060520298044724e-02
   2.1440000000000000e+03   2.8284271247461898e+00                      Inf                      Inf                      Inf                      Inf   1.1084323658254252e-02   9.6681278609461619e-03
   2.6680000000000000e+03   2.8284271247461898e+00                      Inf                      Inf                      Inf                      Inf   8.8378871903355746e-03   7.9598955201620643e-03
   3.6040000000000000e+03   5.6568542494923797e+00                      Inf                      Inf                      Inf                      Inf   7.3155014965563753e-03   6.4888320266014111e-03
   4.4560000000000000e+03   5.6568542494923797e+00                      Inf                      Inf                      Inf                      Inf   5.7728188122327289e-03   5.1502089717731627e-03
   5.5000000000000000e+03   5.6568542494923797e+00                      Inf                      Inf                      Inf                      Inf   4.7775171841083113e-03   4.2276429577163075e-03
   7.0280000000000000e+03   5.6568542494923797e+00                      Inf                      Inf                      Inf                      Inf   3.9236567415101105e-03   3.4905817495464414e-03
   8.9040000000000000e+03   5.6568542494923797e+00                      Inf                      Inf                      Inf                      Inf   3.1619433382232224e-03   2.7930342367230797e-03
   1.0520000000000000e+04   5.6568542494923797e+00                      Inf                      Inf                      Inf                      Inf   2.5942260993471844e-03   2.2750822987596458e-03
   1.2764000000000000e+04   5.6568542494923797e+00                      Inf                      Inf                      Inf                      Inf   2.1285818494612050e-03   1.8779815294061378e-03
   1.5940000000000000e+04   5.6568542494923797e+00                      Inf                      Inf                      Inf                      Inf   1.7727280132718050e-03   1.5733385209995366e-03
   1.9196000000000000e+04   5.6568542494923797e+00                      Inf                      Inf                      Inf                      Inf   1.4734922920531090e-03   1.2987480557937554e-03
   2.4024000000000000e+04   5.6568542494923797e+00                      Inf                      Inf                      Inf                      Inf   1.2668640158449563e-03   1.1298567942328477e-03
   2.8476000000000000e+04   5.6568542494923797e+00                      Inf                      Inf                      Inf                      Inf   1.0365801686428688e-03   9.1830484025059837e-04
   3.5432000000000000e+04   1.1313708498984759e+01                      Inf                      Inf                      Inf                      Inf   8.4726596496623500e-04   7.4866371428906214e-04
   4.3468000000000000e+04   1.1313708498984759e+01                      Inf                      Inf                      Inf                      Inf   6.9486153617097542e-04   6.1374449847216822e-04
   5.2824000000000000e+04   1.1313708498984759e+01                      Inf                      Inf                      Inf                      Inf   5.7587964283184532e-04   5.0740011174961344e-04
   6.4152000000000000e+04   1.1313708498984759e+01                      Inf                      Inf                      Inf                      Inf   4.8011779094721962e-04   4.2307568289760047e-04
   7.7976000000000000e+04   1.1313708498984759e+01                      Inf                      Inf                      Inf                      Inf   4.0170508797550560e-04   3.4559915612318093e-04
   9.6364000000000000e+04   1.1313708498984759e+01                      Inf                      Inf                      Inf                      Inf   1.7378071955981222e-03   1.7287984006784528e-03
   9.7604000000000000e+04   1.1313708498984759e+01                      Inf                      Inf                      Inf                      Inf   3.3563818845584850e-04   2.9456158985147847e-04
   1.1473600000000000e+05   1.1313708498984759e+01                      Inf                      Inf                      Inf                      Inf   2.8640772734787407e-04   2.4625141361349942e-04
    }\pucciadapt
	\centering
	    \begin{minipage}{.66\textwidth}
		\begin{tikzpicture}
		\begin{loglogaxis}[legend pos=south west,legend cell align=left,
			legend style={fill=none,draw=none},
			ymin=1e-4,ymax=1e0,
			xmin=10,xmax=3e5,
			ytick={1e-4,1e-3,1e-2,1e-1,1e0,1e1},
			xtick={1e1,1e2,1e3,1e4,1e5,1e6,1e7}
			]
			\pgfplotsset{
				cycle list={%
					{black, mark=diamond*, mark size=1.5pt},
					{cyan, mark=*, mark size=1.5pt},
					{orange, mark=square*, mark size=1.5pt},
					{blue, mark=triangle*, mark size=1.5pt},
					{black, mark=none, mark size=1.5pt},
					{black, dashed, mark=none, mark size=1.5pt},
					{black, dashed, mark=diamond*, mark size=1.5pt},
					{cyan, dashed, mark=*, mark size=1.5pt},
					{orange, dashed, mark=square*, mark size=1.5pt},
					{blue, dashed, mark=triangle*, mark size=1.5pt},
				},
				legend style={
					at={(0.05,0.05)}, anchor=south west},
				font=\sffamily\scriptsize,
				xlabel=ndof,ylabel=
			}
			\addplot+ table[x=ndof,y=energy]{\pucciadapt};
			\addlegendentry{\tiny $\Phi_j^{1/2}$ adaptive}
			
			\addplot+ table[x=ndof,y=energy]{\pucciunif};
			\addlegendentry{\tiny $\Phi_j^{1/2}$ uniform}
			
			\addplot+ [dotted,mark=none,black] coordinates{(1e0,1.3e0) (1e12,1.3e-3)};
			\addlegendentry{\tiny rate $\mathsf{ndof}^{-{1/4}}$}
			\addplot+ [dashed,mark=none,black] coordinates{(1e0,1e2) (1e10,1e-8)};
			\addlegendentry{\tiny rate $\mathsf{ndof}^{-1}$}
			
		\end{loglogaxis}
	\end{tikzpicture}
	\end{minipage}
	\hfil
	\begin{minipage}{.33\textwidth}
	\includegraphics[width=\textwidth]{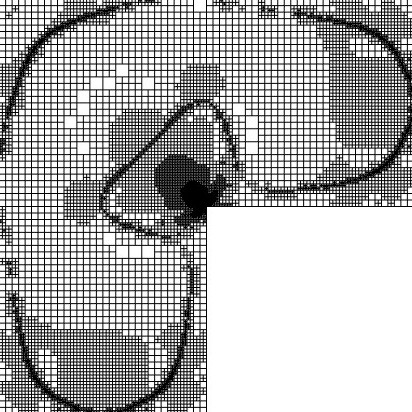}
	\end{minipage}
	\caption{Two-dimensional Pucci experiment from \Cref{ss:numPucci2d}.
	 Convergence history (left) and adaptive mesh (right).}
	\label{fig:pucci_conv}
\end{figure}

\Cref{fig:pucci_conv} displays the suboptimal convergence rate 1/4 due to the expected singularity at the re-entrant corner. The adaptive algorithm refines towards the corner singularity as shown in the mesh
from \Cref{fig:pucci_conv} and almost recovers the optimal convergence rate of 1.

\subsection{Conclusions}
This paper proposes convergent minimal residual methods for fully nonlinear PDE with strong solutions based on the Alexandrov--Bakelman--Pucci maximum principle.
The main challenge is the non-convexity and non-smoothness of the resulting finite dimensional minimization problem. Nevertheless, the proposed policy iteration performs well and the residual serves as an a~posteriori indicator for convergence of the errors. The adaptive algorithm motivated by this residual performs better than uniform mesh refinements for singular solutions.

\bibliographystyle{amsplain}
\bibliography{references.bib}

\providecommand{\bysame}{\leavevmode\hbox to3em{\hrulefill}\thinspace}
\providecommand{\MR}{\relax\ifhmode\unskip\space\fi MR }
\providecommand{\MRhref}[2]{%
  \href{http://www.ams.org/mathscinet-getitem?mr=#1}{#2}
}
\providecommand{\href}[2]{#2}
\begin{thebibliography}{10}

\bibitem{BarlesSouganidis1991}
G.~Barles and P.~E. Souganidis, \emph{Convergence of approximation schemes for
  fully nonlinear second order equations}, Asymptotic Anal. \textbf{4} (1991),
  no.~3, 271--283. \MR{1115933}

\bibitem{BrennerScott2008}
S.~C. Brenner and L.~R. Scott, \emph{The mathematical theory of finite element
  methods}, third ed., Texts in Applied Mathematics, vol.~15, Springer, New
  York, 2008.

\bibitem{BrennerGudiSung2010}
Susanne~C. Brenner, Thirupathi Gudi, and Li-yeng Sung, \emph{An a posteriori
  error estimator for a quadratic {$C^0$}-interior penalty method for the
  biharmonic problem}, IMA J. Numer. Anal. \textbf{30} (2010), no.~3, 777--798.
  \MR{2670114}

\bibitem{CaffarelliCrandallKocanSwiech1996}
L.~Caffarelli, M.~G. Crandall, M.~Kocan, and A.~\'Swi\k{e}ch, \emph{On
  viscosity solutions of fully nonlinear equations with measurable
  ingredients}, Comm. Pure Appl. Math. \textbf{49} (1996), no.~4, 365--397.
  \MR{1376656}

\bibitem{CloughTocher1965}
R.~W. Clough and J.~L. Tocher, \emph{Finite element stiffness matrices for
  analysis of plates in bending}, Proceedings of the Conference on Matrix
  Methods in Structural Mechanics (1965), 515--545.

\bibitem{CrandallIshiiLions1992}
Michael~G. Crandall, Hitoshi Ishii, and Pierre-Louis Lions, \emph{User's guide
  to viscosity solutions of second order partial differential equations}, Bull.
  Amer. Math. Soc. (N.S.) \textbf{27} (1992), no.~1, 1--67. \MR{1118699}

\bibitem{Doktor1976}
Pavel Doktor, \emph{Approximation of domains with {L}ipschitzian boundary},
  \v{C}asopis P\v{e}st. Mat. \textbf{101} (1976), no.~3, 237--255. \MR{0461122}

\bibitem{Doerfler1996}
W.~D{\"o}rfler, \emph{A convergent adaptive algorithm for {P}oisson's
  equation}, SIAM J. Numer. Anal. \textbf{33} (1996), no.~3, 1106--1124.

\bibitem{DouglasDupontPercellScott1979}
Jim Douglas, Jr., Todd Dupont, Peter Percell, and Ridgway Scott, \emph{A family
  of {$C^{1}$} finite elements with optimal approximation properties for
  various {G}alerkin methods for 2nd and 4th order problems}, RAIRO Anal.
  Num\'{e}r. \textbf{13} (1979), no.~3, 227--255. \MR{543934}

\bibitem{ErnGuermond2021}
Alexandre Ern and Jean-Luc Guermond, \emph{Finite elements
  {I}---{A}pproximation and interpolation}, Texts in Applied Mathematics,
  vol.~72, Springer, Cham, 2021. \MR{4242224}

\bibitem{Evans1982}
Lawrence~C. Evans, \emph{Classical solutions of fully nonlinear, convex,
  second-order elliptic equations}, Comm. Pure Appl. Math. \textbf{35} (1982),
  no.~3, 333--363. \MR{649348}

\bibitem{Figalli2017}
Alessio Figalli, \emph{The {M}onge-{A}mp\`ere equation and its applications},
  Zurich Lectures in Advanced Mathematics, European Mathematical Society (EMS),
  Z\"{u}rich, 2017. \MR{3617963}

\bibitem{Gallistl2015}
D.~Gallistl, \emph{{M}orley finite element method for the eigenvalues of the
  biharmonic operator}, IMA J. Numer. Anal. \textbf{35} (2015), no.~4,
  1779--1811.

\bibitem{GallistlTian2024}
D.~Gallistl and S.~Tian, \emph{A posteriori error estimates for nonconforming
  discretizations of singularly perturbed biharmonic operators}, SMAI J.
  Comput. Math. \textbf{10} (2024), 355--372.

\bibitem{GallistlTran2023}
D.~Gallistl and N.~T. Tran, \emph{Convergence of a regularized finite element
  discretization of the two-dimensional {M}onge--{A}mp\`ere equation}, Math.
  Comp. \textbf{92} (2023), no.~342, 1467--1490.

\bibitem{GallistlTran2024}
\bysame, \emph{Stability and guaranteed error control of approximations to the
  {M}onge--{A}mp\`ere equation}, Numer. Math. \textbf{156} (2024), no.~1,
  107--131.

\bibitem{GallistlSueli2019}
Dietmar Gallistl and Endre S\"{u}li, \emph{Mixed finite element approximation
  of the {H}amilton-{J}acobi-{B}ellman equation with {C}ordes coefficients},
  SIAM J. Numer. Anal. \textbf{57} (2019), no.~2, 592--614. \MR{3924618}

\bibitem{GeorgoulisHoustonVirtanen2011}
Emmanuil~H. Georgoulis, Paul Houston, and Juha Virtanen, \emph{An {\it a
  posteriori} error indicator for discontinuous {G}alerkin approximations of
  fourth-order elliptic problems}, IMA J. Numer. Anal. \textbf{31} (2011),
  no.~1, 281--298. \MR{2755946}

\bibitem{GilbargTrudinger2001}
David Gilbarg and Neil~S. Trudinger, \emph{Elliptic partial differential
  equations of second order}, Classics in Mathematics, Springer-Verlag, Berlin,
  2001, Reprint of the 1998 edition. \MR{1814364}

\bibitem{Grisvard2011}
Pierre Grisvard, \emph{Elliptic problems in nonsmooth domains}, Classics in
  Applied Mathematics, vol.~69, SIAM, Philadelphia, PA, 2011. \MR{3396210}

\bibitem{GuzmanLischkeNeilan2022}
Johnny Guzm{\'a}n, Anna Lischke, and Michael Neilan, \emph{Exact sequences on
  {W}orsey--{F}arin splits}, Mathematics of Computation \textbf{91} (2022),
  no.~338, 2571--2608.

\bibitem{Jensen2017}
Max Jensen, \emph{$l^2(h_\gamma^1)$ finite element convergence for degenerate
  isotropic {H}amilton-{J}acobi-{B}ellman equations}, IMA J. Numer. Anal.
  \textbf{37} (2017), no.~3, 1300--1316.

\bibitem{JensenSmears2013}
Max Jensen and Iain Smears, \emph{On the convergence of finite element methods
  for {H}amilton-{J}acobi-{B}ellman equations}, SIAM J. Numer. Anal.
  \textbf{51} (2013), no.~1, 137--162.

\bibitem{KaweckiSmears2021}
Ellya~L. Kawecki and Iain Smears, \emph{Unified analysis of discontinuous
  {G}alerkin and $c^0$-interior penalty finite element methods for
  {H}amilton-{J}acobi-{B}ellman and {I}saacs equations}, ESAIM Math. Model.
  Numer. Anal. \textbf{55} (2021), no.~2, 449--478.

\bibitem{KaweckiSmears2022}
\bysame, \emph{Convergence of adaptive discontinuous {G}alerkin and
  {$C^0$}-interior penalty finite element methods for
  {H}amilton-{J}acobi-{B}ellman and {I}saacs equations}, Found. Comput. Math.
  \textbf{22} (2022), no.~2, 315--364. \MR{4407745}

\bibitem{KoikeSwiech2009}
Shigeaki Koike and Andrzej \'{S}wi\c{e}ch, \emph{Weak {H}arnack inequality for
  fully nonlinear uniformly elliptic {PDE} with unbounded ingredients}, J.
  Math. Soc. Japan \textbf{61} (2009), no.~3, 723--755. \MR{2552914}

\bibitem{Krylov1983}
N.~V. Krylov, \emph{Boundedly inhomogeneous elliptic and parabolic equations in
  a domain}, Izv. Akad. Nauk SSSR Ser. Mat. \textbf{47} (1983), no.~1, 75--108.
  \MR{688919}

\bibitem{McShane1934}
E.~J. McShane, \emph{Extension of range of functions}, Bull. Amer. Math. Soc.
  \textbf{40} (1934), no.~12, 837--842. \MR{1562984}

\bibitem{NeilanWu2019}
Michael Neilan and Mohan Wu, \emph{Discrete {M}iranda-{T}alenti estimates and
  applications to linear and nonlinear {PDE}s}, J. Comput. Appl. Math.
  \textbf{356} (2019), 358--376.

\bibitem{Safonov1988}
M.~V. Safonov, \emph{Classical solution of second-order nonlinear elliptic
  equations}, Izv. Akad. Nauk SSSR Ser. Mat. \textbf{52} (1988), no.~6,
  1272--1287, 1328. \MR{984219}

\bibitem{SmearsSueli2013}
Iain Smears and Endre S\"{u}li, \emph{Discontinuous {G}alerkin finite element
  approximation of nondivergence form elliptic equations with {C}ord\`es
  coefficients}, SIAM J. Numer. Anal. \textbf{51} (2013), no.~4, 2088--2106.
  \MR{3077903}

\bibitem{SmearsSueli2014}
\bysame, \emph{Discontinuous {G}alerkin finite element approximation of
  {H}amilton-{J}acobi-{B}ellman equations with {C}ordes coefficients}, SIAM J.
  Numer. Anal. \textbf{52} (2014), no.~2, 993--1016. \MR{3196952}

\bibitem{Tran2024}
Ngoc~Tien Tran, \emph{Finite element approximation for uniformly elliptic
  linear {PDE} of second order in nondivergence form}, Math. Comp. \textbf{94}
  (2025), no.~353, 1043--1064.

\bibitem{WorseyFarin1987}
A.~J. Worsey and G.~Farin, \emph{An {$n$}-dimensional {C}lough-{T}ocher
  interpolant}, Constr. Approx. \textbf{3} (1987), no.~2, 99--110. \MR{889547}

\end{thebibliography}

\end{document}